\documentclass[12pt]{amsart}
\usepackage{amsfonts,amssymb ,amsthm,latexsym}

\newtheorem{Thm}{Theorem}[section] 
\newtheorem{Cor}[Thm]{Corollary} \newtheorem{Lem}[Thm]{Lemma}
\newtheorem{Prop}[Thm]{Proposition} \newtheorem{Def}[Thm]{Definition}

\newcommand{\no}{\nonumber}  \newcommand{\bea}{\begin{eqnarray}}
\newcommand{\eea}{\end{eqnarray}}

       \def \Bbb{\mathbb} 
\begin{document}

\title{Bracket products for Weyl-Heisenberg frames}
\author{Peter G. Casazza and M. C. Lammers }
 \address{Department of
Mathematics \\ University of Missouri-Columbia \\ Columbia, MO 65211\\
and Department of Mathematics\\ The University of South Carolina\\
Columbia, SC 29208 }

\thanks{The first author was supported by NSF DMS 970618}

\email{pete@math.missouri.edu;lammers@math.sc.edu}



\begin{abstract}
We provide a detailed development
 of a function valued inner product known  as  the bracket
product and used effectively by de Boor, Devore,
Ron and Shen to study translation invariant systems.  We develop
a version of the bracket product
specifically geared to Weyl-Heisenberg frames.    This bracket
product has all the properties of a standard inner product
including Bessel's
inequality, a Riesz Representation Theorem, and a  Gram-Schmidt
process which
turns a sequence of functions $(g_{n})$ into a sequence  $(e_{n})$  
with the
property that  $(E_{mb}e_{n})_{m,n\in \Bbb Z}$ is orthonormal in
$L^{2}(\Bbb
R)$.  Armed with this inner product, we obtain several results
concerning
Weyl-Heisenberg frames.  First we see that fiberization in this setting
takes on a particularly simple form and we use it to obtain a
compressed
representation of the
frame operator.   Next, we write down explicitly all those
functions  $g\in
L^{2}(\Bbb R)$ and $ab=1$ so that the family $(E_{mb}T_{na}g)$ is
complete in $L^{2}(\Bbb R)$.  One consequence of this is that for
functions $g$ supported on a half-line $[{\alpha},\infty)$
(in particular, for compactly supported $g$), $(g,1,1)$
is complete if and only if $\text{sup}_{0\le t< a}|g(t-n)|\not= 0$  
a.e.
Finally, we give a direct proof of a result hidden in the literature
 by proving:  For any $g\in   L^{2}(\Bbb R)$, $A\le
\sum_{n}|g(t-na)|^{2}\le B$ is equivalent to $(E_{m/a}g)$ being a
Riesz basic sequence.
\end{abstract}

\maketitle



\section{Introduction} \label{intro}

While working on some deep questions in non-harmonic Fourier series,
 Duffin and Schaeffer \cite{DS} introduced the notion of a frame for
 Hilbert spaces.  Outside of this area, this idea seems to have been
 lost until Daubechies, Grossman and Meyer \cite{DGM} brought
 attention to it in 1986.  Duffin and Schaeffer's definition was an
 abstraction of a concept introduced by Gabor \cite{G} in 1946 for
 doing signal analysis.  Today the frames introduced by Gabor are
 called {\bf Gabor frames} or {\bf Weyl-Heisenberg frames}.
 Along with wavelets, Weyl-Heisenberg frames are still the backbone
 of modern day signal
 processing as well as a host of related topics.\\

In the study  of shift invariant systems and frames several authors,
including de Boor, DeVore, Ron and Shen \cite{BDR1,BDR2,RS1,RS2},
have made extensive use of the so called {\bf bracket product}
\[ [f,g](x)=\sum_{\beta\in 2 \pi^d}f(x+\beta)\overline{g(x+\beta)}.\]
One may view this bracket product as a pointwise inner product and we
will refer to it as such throughout the paper.  In what follows we
give a more thorough development of the bracket product itself and
its application to univariate principal Weyl-Heisenberg systems. We
hope that our development of the bracket product will aid in applying
it to Weyl-Heisenberg systems as well as other areas where
shift-invariance is of  importance.  Because we would like to be
able to change the shift parameter  from $2\pi$ to arbitrary $a\in
\mathbb R^+$  we will refer to this bracket
product as the {\bf a-inner product}.\\

Let us briefly discuss the organization of the paper.  In Section 2 we
review the notation and terminology, as well as the basic results of
Weyl-Heisenberg frames.  In Section 3 we ever so  slightly alter
 the definition
of bracket product  to get  the a-inner product
and develop its basic properties.  In section 4 we discuss
orthogonality with respect to the $a$-inner product and develop such
notions as orthonormal sequences, orthonormal bases and a
Bessel's inequality all with respect to the a-inner product..
In Section 5 we study a-factorable operators.  These are the natural
bounded linear operators related to the a-inner product.  We will
prove that the a-inner product has a Riesz Representation Theorem for
a-factorable   operators.  In Section 6 we will relate our a-inner
product directly to   Weyl-Heisenberg frames.  We will see that this
gives a  representation for the frame operator for a Weyl-Heisenberg
frame $(g,a,b)$ in terms of the $1/b$-inner product. This
representation can be viewed as a simple  form of fiberization
technique developed by  Ron and Shen \cite{RS1,RS2}.
 In Section 7 we use these ideas to prove two   theorems concerning
Weyl-Heisenberg frames.  The first is
``half'' of a   result proved independently by Daubechies, H. Landau,
Z. Landau \cite{DLL} ; Janssen \cite{J95}; and by Ron and Shen
\cite{RS2}.  The second is a complete listing of all   functions $g\in
L^{2}(\Bbb R)$ and $ab=1$ so that the Weyl-Heisenberg system   is
complete.  A surprising consequence of this is that for a function
supported on a half line, the minimal necessary condition for
completeness $\text{sup}_{n}|g(t-na)|\not= 0$ a.e. becomes
sufficient.  In Section 8 we see that the a-inner product gives a
natural   definition for an a-frame, and that these frames are a
natural generalization of   regular frames.  In particular, we show 
that $(g,a,b)$ is a WH-frame iff the trnslates of g, $(g,a)$, forms a 
(1/b)-frame.  We will also look at a-Riesz   bases and their
relationship to Riesz bases for a Hilbert space.  Finally, in
Section 9 we show that the Gram-Schmidt orthogonalization procedure
works   exactly as expected to produce a-orthonormal sequences with
the proper spans.\\

The authors would like to thank A.J.E.M. Janssen for his helpful
comments.  In particular we would like to thank him for pointing out
the connection between what we refer to as the compression of the
frame operator and fiberization.  Also, we would like to thank
R. DeVore and A. Ron for useful discussions concerning the material
in this paper.



\section{Preliminaries}\label{Prelim}
\setcounter{equation}{0}

We use $\Bbb N, \Bbb Z, \Bbb R, \Bbb C$ to denote the natural numbers,
integers, real numbers and complex numbers, respectively.  A {\bf
scalar} is an element of $\Bbb R$ or $\Bbb C$.  Integration is always
with respect to Lebesgue measure.  $L^{2}(\Bbb R)$ will denote the
complex Hilbert space of square integrable functions mapping   $\Bbb
R$ into $\Bbb C$.  A bounded unconditional basis for a Hilbert space
$H$ is   called a {\bf Riesz basis}.  That is, $(f_{n})$ is a Riesz
basis for $H$ if   and only if there is an orthonormal basis $(e_{n})$
for $H$ and an  operator  $T:H\rightarrow H$ defined by $T(e_{n}) =
f_{n}$, for all $n$.  We call $(f_{n})$ a {\bf Riesz basic sequence}
if it is a Riesz basis for its closed linear span.  For $E\subset H$,
we write span E for the {\bf   closed linear span of E}.\\

In 1952, Duffin and Schaeffer \cite{DS} were working on some deep
problems in non-harmonic Fourier series.  This led them to define

\begin{Def}\label{Frame}  A sequence $(f_{n})_{n\in \Bbb Z}$ of elements
of a Hilbert space $H$ is called a {\bf frame} if there are constants
$A,B>0$ such that
\begin{equation}\label{FrameEqn}
A\|f\|^{2}\le \sum_{n\in \Bbb Z}|<f,f_{n}>|^{2}\le B\|f\|^{2},\ \
\text{for all}\ \ f\in H.
\end{equation}
\end{Def}

The numbers $A,B$ are called the {\bf lower} and {\bf upper frame
bounds} respectively.    The largest number $A>0$ and smallest number
$B>0$ satisfying the frame   inequalities for all $f\in H$ are called
the {\bf optimal frame bounds}.  The frame is a {\bf tight frame} if
$A=B$ and a {\bf normalized tight frame} if $A=B=1$.  A frame is {\bf
exact} if it ceases to be a frame when any one of its elements is
removed.  It is known that a frame is exact if and only if it is a
Riesz basis.  A non-exact frame is called {\bf over-complete}  in the
sense that at   least one vector can be removed from the frame and the
remaining set of   vectors will still form a frame for $H$ (but
perhaps with different frame bounds).  If $f_{n}\in H$, for all $n\in
\Bbb Z$, we call $(f_{n})_{n\in \Bbb Z}$ a {\bf frame sequence} if it
is a frame for its closed linear span in $H$.\\

We will consider frames from the operator theoretic point of view.  To
formulate this approach, let $(e_{n})$ be an orthonormal basis for an
infinite dimensional Hilbert space $H$ and let $f_{n}\in H$, for all
$n\in \Bbb Z$.  We call the operator $T:H\rightarrow H$ given by
$Te_{n}=f_{n}$ the {\bf preframe operator} associated with $(f_{n})$.
Now, for each $f\in H$ and $n\in \Bbb Z$ we have $<T^{*}f,e_{n}> =
<f,Te_{n}>=<f,f_{n}>$.  Thus
\begin{equation} \label{ET}
T^{*}f = \sum_{n}<f,f_{n}>e_{n},\ \ \text{for all} \ \ f\in H.
\end{equation}

By (\ref{ET})
$$
\|T^{*}f\|^{2} = \sum_{n}|<f,f_{n}>|^{2},\ \ \text{for all}\ \ f\in H.
$$
It follows that the preframe operator is bounded if and only if
$(f_{n})$ has a finite upper frame bound $B$.  Comparing this to
Definition \ref{Frame} we have

\begin{Thm}\label{T1}
Let $H$ be a Hilbert space with an orthonormal basis $(e_{n})$. Also
 let $(f_{n})$ be a sequence of elements of $H$ and let $Te_{n}=f_{n}$
 be the preframe operator.  The following are equivalent:

(1)  $(f_{n})$ is a frame for $H$.

(2)  The operator $T$ is bounded, linear and onto.

(3)  The operator $T^{*}$ is an (possibly into) isomorphism called
the {\bf frame transform}.

Moreover, $(f_{n})$ is a normalized tight frame if and only if the
preframe operator is a quotient map (i.e. a co-isometry).
\end{Thm}

The dimension of the kernel of T is called the {\bf excess} of the
frame.  It follows that $S=TT^{*}$ is an invertible operator on   $H$,
called the {\bf frame operator}.  Moreover, we have
$$
Sf = TT^{*}f = T(\sum_{n}<f,f_{n}>e_{n}) = \sum_{n}<f,f_{n}>Te_{n}
=\sum_{n}<f,f_{n}>f_{n}.
$$

A direct calculation now yields
$$
<Sf,f> = \sum_{n}|<f,f_{n}>|^{2}.
$$
Therefore, the {\bf frame operator is a positive, self-adjoint
invertible operator} on $H$.  Also, the frame inequalities
(\ref{Frame}) yield that $(f_{n})$ is a frame with frame bounds
$A,B>0$ if and only if $A\cdot I\le S\le B\cdot I$.  Hence, $(f_{n})$
is a normalized tight frame if and only if $S=I$.  Also, a direct
calculation yields \bea\label{E1} f = SS^{-1}f &=&
\sum_{n}<S^{-1}f,f_{n}>f_{n} \\ &=& \sum_{n}<f,S^{-1}f_{n}>f_{n}\no \\
&=&\sum_{n}<f,S^{-1/2}f_{n}>S^{-1/2}f_{n}.\no \eea

We call $(<S^{-1}f,f_{n}>)$ the {\bf frame coefficients} for $f$.  One
interpretation of equation (\ref{E1}) is that $(S^{-1/2}f_{n})$ is a
normalized tight frame.

\begin{Thm}\label{Tight}
Every frame $(f_{n})$ (with frame operator $S$) is equivalent to the
normalized tight frame $(S^{-1/2}f_{n})$.
\end{Thm}

We will work here with a particular class of frames called
Weyl-Heisenberg frames.  To formulate these frames, we first need some
notation.  For a function $f$ on $\Bbb R$ we define the operators:

\[
\begin{array}{lll}
\text{Translation:} & T_{a}f(x) = f(x-a), & a\in \Bbb R \\
\text{Modulation:}  & E_{a}f(x) = e^{2{\pi}iax}f(x), & a\in \Bbb R \\
\text{Dilation:} & D_{a}f(x) = |a|^{-1/2}f(x/a), & a\in \Bbb {R} -
\{0\}
\end{array}
\]

We also use the symbol $E_{a}$ to denote the {\bf exponential
function} $E_{a}(x) = e^{2{\pi}iax}$.  Each of the operators $T_{a},
E_{a}, D_{a}$ are unitary operators on  $L^{2}(\Bbb R )$ and they
satisfy:

\[
\begin{array}{l}
T_{a}E_{b}f(x) = e^{2{\pi}ib(x-a)}f(x-a); \\ E_{b}T_{a}f(x) =
e^{2{\pi}ibx}f(x-a); \\ D_{a}T_{b}f(x) = |a|^{-1/2}f(\frac{x}{a}-b);
\\ T_{b}D_{a}f(x) = |a|^{-1/2}f(\frac{x-b}{a}); \\ E_{b}D_{a}f(x) =
e^{2{\pi}ibx}|a|^{-1/2}f(\frac{x}{a}); \\ D_{a}E_{b}f(x) =
e^{2{\pi}ibx/a}|a|^{-1/2}f(\frac{x}{a}).
\end{array}
\]

In 1946 Gabor \cite{G} formulated a fundamental approach to signal
decomposition in terms of elementary signals.  This method resulted in
{\bf Gabor frames} or as they are often called today {\bf
Weyl-Heisenberg frames}.

\begin{Def}\label{WHS}
If $a,b\in \Bbb R$ and $g\in L^{2}(\Bbb R )$ we call
$(E_{mb}T_{na}g)_{m,n\in \Bbb Z}$ a {\bf Weyl-Heisenberg system} ({\bf
WH-system} for short) and denote it by $(g,a,b)$.  We denote by $(g,a)$
the family $(T_{na}g)_{n\in \mathbb Z}$.
We call $g$ the
{\bf window function}.
\end{Def}

If the WH-system $(g,a,b)$ forms a frame for $L^{2}(\Bbb R)$, we
call this a {\bf Weyl-Heisenberg frame} ({\bf WH-frame} for short).
The   numbers $a,b$ are the {\bf frame parameters} with $a$ being the
{\bf shift parameter}   and $b$ being the {\bf modulation parameter.}
We will be interested in   when there are finite upper frame bounds
for a WH-system.  We call this class of functions the {\bf preframe
functions} and denote this class by {\bf PF}.  It is easily checked
that

\begin{Prop}\label{P6}
The following are equivalent:

(1) $g\in$ {\bf PF}.

(2)  The operator
$$
Sf = \sum_{n}<f,E_{mb}T_{na}g>E_{mb}T_{na}g,
$$
is a well defined bounded linear operator on $L^{2}(\Bbb R)$.
\end{Prop}

We will need the WH-frame identity due to Daubechies \cite{D}.  To
simplify the notation a little we introduce the following auxiliary
functions defined for a $g\in L^{2}(\Bbb R)$ and all $k\in \Bbb Z$ by

$$
G_{k}(t) = \sum_{n\in \Bbb Z} g(t-na)\overline{ g(t-na-k/b)}.
$$

In particular,
$$
G_{0}(t) = \sum_{n\in \Bbb Z}|g(t-na)|^{2}.
$$

\begin{Thm}\label{WHFI}({\bf WH-Frame Identity.})
If $\sum_{n}|g(t-na)|^2 \le B$ a.e. and $f\in L^{2}(\Bbb R)$ is
bounded and compactly supported, then
$$
\sum_{n,m\in \Bbb Z}|<f,E_{mb}T_{na}g>|^{2} = F_{1}(f)+F_{2}(f),
$$
where
$$
F_{1}(f) = b^{-1}\int_{R}|f(t)|^{2}G_{0}(t)\  dt,
$$
and \bea F_{2}(f)  &=& b^{-1}\sum_{k\not=
0}\int_{R}\overline{f(t)}f(t-k/b)G_{k}(t)\ dt \no \\
&=&b^{-1}\sum_{k\ge
1}2\text{Re}\int_{R}\overline{f(t)}f(t-k/b)G_{k}(t)\ dt.  \no \eea
\end{Thm}

There are many restrictions on the $g,a,b$ in order that $(g,a,b)$
form a WH-frame.  We will make use of a few of them here.  The first
is a simple application of the WH-frame Identity.  That is, if we put
functions supported on $[0,1/b]$ into this identity, then $F_{2}(f) =
0$.  Now the WH-frame Identity combined with the frame condition
quickly yields,

\begin{Thm}\label{AB}
If $(g,a,b)$ is a WH-frame with frame bounds $A,B$ then
$$
A\le bG_{0}(t)\le B,\ \ \text{a.e.}
$$
\end{Thm}

Casazza and Christensen \cite{CC1} noted that we have a similar upper
bound condition with a replaced by 1/b.

\begin{Prop}\label{PP6}
If $(g,a,b)$ is a WH-frame with upper frame bound B than
$$
\sum_{n\in \Bbb Z}|g(t-n/b)|^{2} \le B,\ \ \text{a.e.}.
$$
\end{Prop}

There are also some restrictions on $a,b$ for $(g,a,b)$ to be a frame.

\begin{Prop}\label{P11}
Let $g\in L^{2}(\Bbb R)$ and $a,b\in \Bbb R$.

(1)  If $(E_{mb}T_{na}g)$ is complete, then $ab\le 1$.

(2)  If $(g,a,b)$ is a WH-frame and

\ \ \ \ (i)  $ab<1$ then $(g,a,b)$ is over-complete.

\ \ \ \ (ii)  $ab=1$ then $(g,a,b)$ is a Riesz basis.
\end{Prop}

Part (1) of Proposition \ref{P11} has a complicated history (see
\cite{D} for a discussion) which derives from the work of Rieffel
\cite{Rieffel}.  Today, there is a simpler   proof using Beurling
density due to Ramanathan and Steger \cite{Steger}.  Moreover, the
results of Ramanathan and Steger \cite{Steger} combined with an
important example of Benedetto, Heil and Walnut \cite{BHW} shows that
the   form of the lattice in the Rieffel result \cite{Rieffel} is
quite important to the  conclusion.  There are many derivations
available for (2) \cite{CCJ1,D,D92,HL,J94,J95}.\\

A recent very important result was proved independently by Daubechies,
H. Landau and Z. Landau \cite{DLL}, Janssen \cite{J95}, and Ron and
Shen \cite{RS2}.

\begin{Thm}\label{dual}
For $g\in L^{2}(\Bbb R)$ and $a,b\in \Bbb R$, the following are
equivalent:

(1)  $(g,a,b)$ is a WH-frame.

(2)  The family $(E_{m/a}T_{n/b}g)_{m,n\in \Bbb Z}$ is a Riesz basic
sequence in $L^{2}(\Bbb R)$.
\end{Thm}

Ron and Shen  attained this result with a technique they call
{\bf Gramian analysis}.  At the heart of this technique is the
Gramian matrix
$\mathcal G$ which is used  to  decompose the pre-frame operator
and its adjoint. The tie in with the bracket product  becomes  clear
when one sees that in the shift-invariant case (i.e. consider only  
$({T_{na}g}$)  this matrix becomes $\mathcal G =[g,g]$.\\

Finally, we will need the classification of tight WH-frames.  Parts of
this are due to various authors.  A direct proof from the definitions
as well as the historical development can be found in \cite{CCJ1}.

\begin{Thm}\label{TightWHF}
Let $g\in L^{2}(\Bbb R)$ and $a,b\in \Bbb R$.  The following are
equivalent:

(1)  $(E_{mb}T_{na}g)$ is a normalized tight Weyl-Heisenberg frame
for $L^{2}(\Bbb R)$.

(2)  We have:

\ \ \ \ (a)  $G_{0}(t) = \sum_{n\in \Bbb Z}|g(t-na)|^{2} = b$ a.e.

\ \ \ \ (b)  For all $k\not= 0$, $G_{k}(t) =
\sum_{n}g(t-na)\overline{g(t-na-k/b)} = 0$ a.e.

(3)  We have $g\perp E_{n/a}T_{m/b}g$, for all $(n,m)\not= (0,0)$
and $\|g\|^{2} = ab$.

(4)  $(E_{n/a}T_{m/b}g)$ is an orthogonal sequence in $L^{2}(\Bbb
R)$ and $\|g\|^{2} = ab$.

(5)  $(E_{mb}T_{na}g)$ is a Weyl-Heisenberg frame for $L^{2}(\Bbb
R)$ with frame operator $S$ and $Sg=g$.

Moreover, when at least one of $(1)-(5)$ holds, $(E_{mb}T_{na}g)$   is
an orthonormal basis for $L^{2}(\Bbb R)$ if and only if $\|g\|=1$.
\end{Thm}

We next recall the {\bf Wiener amalgam space} $W(L^{\infty},L^{1})$
which consists of all functions $g$ so that for some $a>0$ we have,
$$
\|g\|_{W,a} = \sum_{n\in \Bbb Z}\|g\cdot
{\chi}_{[an,a(n+1))}\|_{\infty} = \sum_{n\in \Bbb Z}\|T_{na}\cdot
{\chi}_{[0,a)}\|_{\infty} < \infty.
$$
It is easily checked that $W(L^{\infty},L^{1})$ is a Banach space
with the above norm.  Also, if $\|g\|_{W,a}< \infty$, for one $a>0$,
then   this norm is finite for all $a>0$.



\section{Pointwise Inner Products}\label{PIP}
\setcounter{equation}{0}
A number of the basic results in this section  can be found in various
other papers \cite{BDR1,BDR2,RS1,RS2}.  For the sake of completeness,
and to create a good reference for this inner product we present
them here.
To guarantee that our inner product is well defined, we need
to first check some convergence properties for elements of $L^{2}(\Bbb
R)$.

\begin{Prop}\label{LBD}
For $f,g\in L^{2}(\Bbb R)$ and $a\in \Bbb R$ the series
$$
\sum_{n\in \Bbb Z}f(t-na)\overline{g(t-na)}
$$
converges unconditionally a.e. to a function in $L^{1}[o,a]$.
\end{Prop}

{\it Proof.}  If $f,g\in L^{2}(\Bbb R)$ then $fg\in L^{1}(\Bbb R)$.
Hence, \bea \|fg\|_{L^{1}} &=& \int_{\Bbb R}|f(t)\overline{g(t)}|\ dt
\no \\ & =& \sum_{n\in \Bbb Z}\int_{0}^{a}|f(t-na)\overline{g(t-na)}|\
dt \no \\ &=& \int_{0}^{a}\sum_{n\in \Bbb
Z}|f(t-na)\overline{g(t-na)}|\ dt  < \infty. \no \eea

The last inequality follows by the Monotone Convergence Theorem.  This
yields both the interchange of the integral and the sum and the
existence of $\sum f(t-na)\overline{g(t-na)}$ as a function in
$L^{1}[0,a]$.  \qed \vspace{14pt}

A simple application of the Lebesgue Dominated Convergence Theorem
combined with Proposition \ref{LBD} yields,

\begin{Cor}\label{sums}
For all $f,g\in L^{2}(\Bbb R)$ we have
$$
<f,g> = \int_{0}^{a}\sum_{n\in \Bbb Z}f(t-na)\overline{g(t-na)}\ dt.
$$
\end{Cor}

Now we introduce the pointwise inner product for WH-frames.  We can
also view this as a vector-valued inner product.

\begin{Def}
Fix $a\in \Bbb R$.  For all $f,g\in L^{2}(\Bbb R)$ we define the {\bf
a-pointwise inner product of f and g} (called the {\bf a-inner
product} for short) by
$$
<f,g>_{a}(t) = \sum_{n\in \Bbb Z}f(t-na)\overline{g(t-na)}, \ \
\text{for all}\ \ t\in \Bbb R.
$$
We define the {\bf a-norm of f} by
$$
\|f\|_{a}(t) = \sqrt{<f,f>_{a}}(t).
$$
\end{Def}

We emphasize here that the a-inner product and the a-norm are {\it
functions} on $\Bbb R$ which are clearly a-periodic.  To cut down on
notation,   whenever we have an a-periodic function on $\Bbb R$, we
will also consider it a function on $[0,a]$.   The convergence of
these series is guaranteed by our earlier discussion.  In fact, the
a-inner product $<\cdot , \cdot >_{a}$ is a mapping from $L^{2}(\Bbb
R)\oplus   L^{2}(\Bbb R)$ to the a-periodic functions on $\Bbb R$
whose restriction to   $[0,a]$ lie in $L^{1}[0,a]$. \\

First we show that the a-inner product really is a good generalization
of the standard notion of inner products for a Hilbert space.

\begin{Thm}\label{IP}
Let $f,g,h\in L^{2}(\Bbb R)$, $c,d\in \Bbb C$, and $a,b\in \Bbb R$.
The following properties hold:
\vspace{10pt}

(1)  $<f,g>_{a}$ is a periodic function of period a on $\Bbb R$ with
$<f,g>_{a} \in L^{1}[0,a]$.
\vspace{10pt}

(2)  We have
$$
\|f\|_{L^{2}(\Bbb R)} = \bigg \|\|f\|_{a}(t) \bigg \|_{L^{2}[0,a]}.
$$
\vspace{10pt}

(3)  We have
$$
<f,g> = \int_{0}^{a}<f,g>_{a}(t)\ dt.
$$
\vspace{10pt}

(4)  $<cf+dg,h>_{a} = c<f,h>_{a} + d<g,h>_{a}$.
\vspace{10pt}

(5)  $<f,cg+dh>_{a} = \overline{c}<f,g>_{a} + \overline{d}<f,h>_{a}$.
\vspace{10pt}

(6)  $<f,g>_{a} = \overline{<g,f>_{a}}$.
\vspace{10pt}

(7)  $<fg,h>_{a} = <f,\overline{g}h>_{a}$.
\vspace{10pt}

(8)  If $<f,g>_{a} = 0$ then $<f,g>=0$.
\vspace{10pt}

(9)  $<T_{b}f,T_{b}g>_{a} = T_{b}<f,g>_{a}$.
\vspace{10pt}

(10)  $\|T_{b}g\|_{a}^{2} = T_{b}\|g\|_{a}^{2}$
 \vspace{10pt}

(11)  $<T_{b}f,g>_{a} = T_{b}<f,T_{-b}g>_{a}$.
\vspace{10pt}

(12)  $<f,g>_{a} = \frac{1}{\sqrt{ab}} D_{ab}<D_{\frac{1}{ab}}f,
D_{\frac{1}{ab}}g>_{\frac{1}{b}}$.
\end{Thm}

{\it Proof.}  All the proofs follow directly from the definitions.  We
will give a small sample to show how they proceed.
\vspace{10pt}

(3)  This is just Corollary \ref{sums}.
\vspace{10pt}

(4)  We calculate: \bea  <cf+dg,h>_{a}(t) &=& \sum_{n\in \Bbb
Z}[cf+dg](t-na)\overline{h(t-na)}\no \\ &=&  c\sum_{n\in \Bbb
Z}f(t-na)\overline{h(t-na)} + d\sum_{n\in \Bbb Z}
g(t-na)\overline{h(t-na)} \no \\ &=& c<f,h>_{a}+ d <g,h>_{a}. \no \eea

\vspace{10pt}

(8)  If $<f,g>_{a} = 0$ then by (3),
$$
<f,g> = \int_{0}^{a}<f,g>_{a}(t)\ dt = 0.
$$
\vspace{10pt}

(11) Again, we calculate \bea <T_{b}f,g>_{a} &=& \sum_{n\in \Bbb
Z}f(t-b-na)\overline{g(t-na)} \no \\ &=& T_{b}\sum_{n\in
Z}f(t-na)\overline{g(t-na+b)} = T_{b}<f,T_{-b}g>_{a}.  \no \eea
\vspace{10pt}

(12)  We compute, \bea
<D_{\frac{1}{ab}}f,D_{\frac{1}{ab}}g>_{\frac{1}{b}}(t) &=&
<\sqrt{ab}f(ab\ \cdot ),\sqrt{ab}g(ab\ \cdot )>_{\frac{1}{b}}(t) \no
\\ &=& ab\sum_{n\in \Bbb Z}f(ab(t-n/b))\overline{g(ab(t-n/b))} \no \\
&=& ab\sum_{n\in \Bbb Z}f(abt-na)\overline{g(abt-na)} \no \\ &=&
ab<f,g>_{a}(abt) = \sqrt{ab}D_{\frac{1}{ab}}<f,g>_{a}(t). \no \eea
\qed \vspace{14pt}

Once one sees what is going on, it is not difficult to mimic the
standard proofs for the usual inner product on a Hilbert space to
obtain the following results for the a-inner product.

\begin{Prop}
For all $f,g\in L^{2}(\Bbb R)$ we have,
\vspace{10pt}

(1)  $|<f,g>|_{a}\le \|f\|_{a}\|g\|_{a},$ a.e.
\vspace{10pt}

(2)  $\|f+g\|_{a}^{2} = \|f\|_{a}^{2} + 2Re<f,g>_{a} + \|g\|_{a}^{2}$.
\vspace{10pt}

(3)  $\|f+g\|_{a} \le \|f\|_{a}+\|g\|_{a}$.
\vspace{10pt}

(4)  $\|f+g\|_{a}^{2} + \|f-g\|_{a}^{2} = 2(\|f\|_{a}^{2}
+\|g\|_{a}^{2})$, a.e.
\end{Prop}

Since our a-inner product is an a-periodic function, it enjoys some
special properties related to a-periodic functions.

\begin{Prop} \label{a-per}
Let $f,g\in L^{2}(\Bbb R)$ and let $h\in L^{\infty}(\Bbb R)$ be an
a-periodic function.  Then
$$
<fh,g>_{a} = h<f,g>_{a}\ \ \text{and}\ \ <f,hg>_{a} = \overline{h}
<f,g>_{a}.
$$
In particular, if $h$ satisfies $h(t)\not= 0$ a.e., then $<f,g>_{a} =
0$ if and only if $<fh,g>_{a} = <f,g\overline{h}>_{a} = 0$.
\end{Prop}

{\it Proof.}  We compute \bea <fh,g>_{a}(t)  &=& \sum_{n\in \Bbb
Z}f(t-na)h(t-na)\overline{ g(t-na)}\no \\ &=& \sum_{n\in \Bbb
Z}f(t-na)h(t)\overline{g(t-na)}\no \\  &=&h(t)\sum_{n\in \Bbb
Z}f(t-na)\overline{g(t-na)}  = h(t)<f,g>_{a}(t). \no \eea \qed
\vspace{14pt}

Next we normalize our functions in the a-inner product.  For $f\in
   L^{2}(\Bbb R)$, we define the {\bf a-pointwise normalization of f}
   to be \[ N_{a}(f)(t) = \left \{
\begin{array}{cl}
 \frac{f(t)}{\|f\|_{a}(t)} : & \|f\|_{a}(t)\not= 0 \\ 0 : &
\|f\|_{a}(t) = 0.
\end{array} \right.
\]

We now have

\begin{Prop}
Let $f,g\in L^{2}(\Bbb R)$.

(1)  We have
$$
<N_{a}(f),g>_{a} = \frac{<f,g>_{a}}{\|f\|_{a}},\ \ \text{where}\ \
\|f\|_{a}\not= 0.
$$
In particular, $<f,g>_{a} = 0$ if and only if $<N_{a}(f),g>_{a} = 0$.

(2)  For $f\not= 0$ a.e. we have
$$
<N_{a}(f),N_{a}(f)>_{a}(t) = \sum_{n\in \Bbb Z}|N_{a}(f)(t-na)|^{2}
= 1, a.e.
$$

(3)   we have
$$
\|N_{a}(f)\|^{2}_{L^{2}(\Bbb R)} = {\lambda}(\text{supp}\
\|f\|_{a}|_{[0,a]}) \le a.
$$
where $\lambda$ denotes Lebesgue measure.

(4)  $N_{a}(N_{a}(f)) = N_{a}(f)$.
\end{Prop}

{\it Proof.}  (1)  We compute
$$
<N_{a}(f),g>_{a} = \sum_{n\in \Bbb Z}N_{a}(f)(t-na)\overline{ g(t-na)}
= \sum_{n\in \Bbb Z}\frac{f(t-na)}{\|f\|_{a}(t-na)} \overline{g(t-na)}.
$$
Since our inner product is a-periodic, this equality becomes,
$$
\frac{1}{\|f\|_{a}(t)}\sum_{n\in \Bbb Z}f(t-na)\overline{ g(t-na)} =
\frac{<f,g>_{a}(t)}{\|f\|_{a}(t)},\ \ \text{where}\ \
\|f\|_{a}(t)\not= 0.
$$

(2)  This is two applications of part (1).

(3)  By (2) we have \bea \|N_{a}(f)\|^{2}_{L^{2}(\Bbb R)}  &=&
\int_{\Bbb R}|N_{a}(f)(t)|^{2}\ dt \no \\ &=& \int_{0}^{a}\sum_{n\in
\Bbb Z}|N_{a}(f)(t-na)|^{2}\ dt  = \int_{0}^{a}{\bf 1}_{\text{supp}\
\|f\|_{a}}(t)\ dt\le a. \no \eea

(4)  This is immediate from (2).  \qed \vspace{14pt}


\section{$a$-orthogonality}
The notion of orthogonality with respect to
the a-inner product has been used primarily
 to describe the orthogonal complement in the usual inner product
for shift-invariant spaces.
In this section we explore more thoroughly what it means to be
a-orthogonal and develop such things as a-orthonormal sequences and  a
Bessel inequality for the a-inner product.
This property gives one of the main applications of
the   a-inner product in Weyl-Heisenberg frame theory.  For as we
will see, orthogonality in
this form is very strong.

\begin{Def}
For $f,g\in L^{2}(\Bbb R)$, we say that $f$ and $g$ are {\bf
a-orthogonal}, and write $f{\perp}_{a} g$, if $<f,g>_{a} = 0$.  We
define the {\bf a-orthogonal complement}  of $E\subset L^{2}(\Bbb R)$
by
$$
E^{{\perp}_{a}} = \{g : <f,g>_{a}=0,\ \text{for all}\ f\in E\}.
$$
Similarly, an {\bf a-orthogonal sequence} is a sequence $(f_{n})$
satisfying $f_{n}{\perp}_{a}f_{m}$, for all $n\not= m$.  This is an
{\bf a-orthonormal sequence} if we also have $\|f\|_{a} = 1$, a.e.
where $\|f\|_{a}\not= 0$.  \end{Def}

We now identify an important class of functions for working with the
a-inner product.

\begin{Def}
We say that $g\in L^{2}(\Bbb R)$ is {\bf a-bounded} if there is a
$B>0$ so that
$$
|<g,g>_{a}|\le B,\ \ \text{a.e.}
$$
We let $L^{\infty}_{a}(\Bbb R)$ denote the family of a-bounded
functions.
\end{Def}

We have that $L^{\infty}_{a}(\Bbb R)$ is a non-closed linear subspace
of $L^{\infty}(\Bbb R)$.  To see this, first observe that
$L^{\infty}_{a}(\Bbb R)$ is just the family of functions $g\in
L^{2}(\Bbb R)$ for which $\|g\|_{a}$ is bounded.  So by the properties
we have   developed for $\|\cdot \|_{a}$ we have that
$L^{\infty}_{a}(\Bbb R)$ is a subspace of $L^{\infty}(\Bbb R)$.  Since
$L^{\infty}_{a}(\Bbb R)$ contains all   bounded compactly supported
functions in $L^{2}(\Bbb R)$, and it is easily seen to not equal
$L^{2}(\Bbb R)$, we have that this is a non-closed   subspace.  Note
also that the Wiener amalgam space is a subspace of
$L^{\infty}_{a}(\Bbb R)$.\\

We have not defined orthonormal bases for the a-inner product yet
since, as we will see, this requires a little more care.  First we
need to develop the basic properties of a-orthogonality.

\begin{Prop}\label{PP2}
If $E\subset L^{2}(\Bbb R)$,
$$
E^{{\perp}_{a}} = \cap_{\phi \in L^{\infty}_{a}(\Bbb
R)}({\phi}E)^{\perp} = (\text{span}_{\phi \in L^{\infty}_{a}(\Bbb
R)}{\phi}E)^{\perp}.
$$
\end{Prop}

{\it Proof.}  Let $f\in E^{{\perp}_{a}}$.  For any $g\in E$ and any
a-periodic function $\phi \in L^{\infty}_{a}(\Bbb R)$ we have by
Proposition \ref{a-per}
$$
<f,{\phi}g>_{a} = \overline{\phi}<f,g>_{a} = 0.
$$
Hence, $f{\perp}_{a}{\phi}g$.  That is, $f\in ({\phi}E)^{\perp}$.

Now let $f\in \cap ({\phi}E)^{\perp}$, the intersection being taken
over all bounded a-periodic ${\phi}$.  Let $g\in E$ and define for
$n\in \Bbb N$,
\[
{\phi}_{n}(t) = \left \{
\begin{array}{cl}
 <f,g>_{a}(t) : & |<f,g>_{a}(t)|\le n \\ 0 : & \text{otherwise}.
\end{array} \right.
\]
Note that ${\phi}_{n}$ is a-periodic. Now we compute, \bea 0 =
<f,{\phi}_{n}g>  &=& \int_{\Bbb R}f(t)\overline{{\phi}_{n}(t)}
\overline{g(t)} dt \no \\ &=& \int_{0}^{a}\left ( \sum_{n\in \Bbb
Z}f(t-na) \overline{g(t-na)}\right )  \overline{{\phi}_{n}(t)}\ dt \no
\\ &=& \int_{0}^{a}<f,g>_{a}(t)\overline{{\phi}_{n}(t)}\  dt =
\int_{0}^{a}|{\phi}_{n}(t)|^{2}\ dt . \no \eea Therefore, ${\phi}_{n}
= 0$, for all $n\in \Bbb Z$.  Hence,   $<f,g>_{a} = 0$, and so
$f{\perp}_{a} g$.  That is, $f{\perp}_{a}E$.  \qed \vspace{14pt}

By Theorem \ref{IP} (8), we have that $E^{{\perp}_{a}} \subset
E^{\perp}$.

\begin{Cor}
For $E\subset L^{2}(\Bbb R)$, $E^{{\perp}_{a}}$ is a norm closed
linear subspace of $E^{\perp}$.
\end{Cor}

The  next result which can be found in \cite{BDR1} shows more clearly
what orthoganlity means in this setting .

\begin{Prop}\label{PP1}
For $f,g\in L^{2}(\Bbb R)$, the following are equivalent:

(1)  $f{\perp}_{a} g$.

(2)  $\text{span}_{m\in \Bbb Z}E_{m/a}f \perp \text{span}_{m\in \Bbb
Z}E_{m/a}g$. \end{Prop}

{\it Proof.}  Fix $m\in \Bbb Z$ and compute
$$
<f,E_{m/a}g> = \int_{\Bbb R}f(t)\overline{g(t)}e^{-2{\pi}i(m/a)t}\ dt
=  \int_{0}^{a}<f,g>_{a}(t)e^{-2{\pi}i(m/a)t}\ dt.
$$
It follows that $<f,E_{m/a}g>=0$, for all $m\in \Bbb Z$ if and only if
$<f,g>_{a} = 0$.  A moment's reflection should convince the reader
that this is all we need.  \qed \vspace{14pt}

\begin{Def}
We say that $E\subset L^{2}(\Bbb R)$ is an {\bf a-periodic closed set}
if for any $f\in E$ and any $\phi\in L^{\infty}_{a}(\Bbb R)$ we have
that ${\phi}f\in E$.
\end{Def}

The next result follows immediately from Propositions \ref{PP1} and
\ref{PP2}.

\begin{Cor}
For any $E\subset L^{2}(\Bbb R)$, $E^{{\perp}_{a}}$ is an   a-periodic
closed set. If $E$ is an a-periodic closed set then $E^{\perp} =
E^{{\perp}_{a}}$.
\end{Cor}

Now we observe what orthogonality means for $(E_{m/a}g)$ in terms of
the regular inner product.

\begin{Prop}\label{PP3a}
If $g\in L^{2}(\Bbb R)$ and $\|g\|_{a} = 1$ a.e., then
$(\frac{1}{\sqrt{a}}E_{m/a}g)_{m\in \Bbb Z}$ is an orthonormal
sequence in $L^{2}(\Bbb R)$.
\end{Prop}

{\it Proof.}  For any $n,m\in \Bbb Z$ we have \bea
<E_{n/a}g,E_{m/a}g>  &=& \int_{\Bbb R}|g(t)|^{2}e^{2{\pi}i[(n-m)/a]t}\
dt \no \\ & =& \int_{0}^{a}\|g\|_{a}^{2}(t)e^{2{\pi}i[(n-m)/a]t}\ dt
\no \\ & =& \int_{0}^{a} e^{2{\pi}i[(n-m)/a]t}\ dt =
a{\delta}_{nm}. \no  \eea \qed \vspace{14pt}

\begin{Cor}\label{OR}
If $(g_{n})_{n\in \Bbb N}$ is an a-orthonormal sequence in
$L^{2}(\Bbb R)$, then \\ $(E_{m/a}g_{n})_{n,m\in \Bbb Z}$ is an
orthonormal sequence in $L^{2}(\Bbb R)$.
\end{Cor}

{\it Proof.}  We need that for all $(n,m)\not= ({\ell},k)\in {\Bbb
Z}\times {\Bbb Z}$, $E_{m/a}g_{n}{\perp} E_{{\ell}/a}g_{k}$.  But, if
$m\not= k$, this is Proposition \ref{PP1}, and if $m=k$, this is
Proposition   \ref{PP3a}.  \qed \vspace{14pt}

Corollary \ref{OR} tells us how to define an a-orthonormal basis.

\begin{Def}
Let $g_{n}\in L^{2}(\Bbb R)$.  We call $(g_{n})$ an {\bf a-orthonormal
basis} for $L^{2}(\Bbb R)$ if it is an a-orthonormal sequence and
$$
\overline{\text{span}}\ (E_{m/a}g_{n})_{n,m\in \Bbb Z}= L^{2}(\Bbb R).
$$
\end{Def}

\begin{Prop}
A sequence $(g_{n})$ in $L^{2}(\Bbb R)$ is an a-orthonormal basis if
and only if $(E_{m/a}g_{n})_{n,m\in \Bbb Z}$ is an orthonormal basis
for $L^{2}(\Bbb R)$.
\end{Prop}

We would  like to  capture the important Bessel's Inequality
for a-orthonormal sequences  but before we do so we need to insure that 
$<f,g>_ag$ remains an $L^2(\mathbb R)$ for  functions $g
\in L^{\infty}_{a}(\Bbb R)$.

\begin{Prop}\label{prebessel} If  $g  \in L^{\infty}_{a}(\Bbb R)$ then 
$<f,g>_ag \in L^2(\mathbb R)$ for all $f \in L^2(\mathbb R)$.
\end{Prop}

{\it Proof} First we need to show $<f,g>_a \in L^2([0,a])$.  This follows
 from the Cauchy-Schwarz inequality for the a-inner product.

\bea  \|<f,g>_a\|_{L^{2}[0,a]}^{2} &=& \int_{0}^{a}|<f,g>_{a}(t)|^{2}\ dt
 \no \\ & \le & \int_{0}^{a}<f,f>_{a}(t)<g,g>_{a}(t)\ dt \no \\ &\le &
 B\int_{0}^{a}<f,f>_{a}(t)\ dt  = B\|f\|_{L^{2}(\Bbb R)}^{2}. \no  \eea

Now we can get our results which follows from the Monotone convergence 
theorem and the result above.

\bea \|<f,g>_a g\|_{L^{2}(\mathbb R)}^{2}
 &=& \int_{\mathbb R} |<f,g>_{a}(t)g(t)|^{2}\ dt \no \\
 &=& \sum \int_0^a  |<f,g>_{a}(t)|^{2} |g(t-na|^2 dt \no \\
& \le& \int_0^a  |<f,g>_{a}(t)|^{2} <g,g>_a(t) dt \no \\
& \le & B^2 \|f\|^2_{L^{2}(\Bbb R)} \no \eea

\qed \vspace{14pt}

\begin{Thm}\label{Bessel}
If $(g_{n})_{n\in \Bbb N}$ is an a-orthonormal sequence in
$L^{2}(\Bbb R)$, then for all $f\in L^{2}(\Bbb R)$ we have that

(1) the series of functions $\sum_{n\in \Bbb N}<f,g_{n}>_{a}g_{n}$
converges in $L^{2}(\Bbb R)$.
\vspace{14pt}

(2)  We have ``Bessel's Inequality'',
$$
<f,f>_{a}\ \ \ge \ \ \sum_{n=1}^{\infty}|<f,g_{n}>_{a}|^{2}.
$$
Note that this is an inequality for functions.

Moreover, if $f\in \text{span}\ (E_{m/a}g_{n})_{m,n\in \Bbb Z}$, then
$$
<f,f>_{a} =  \sum_{n=1}^{\infty}|<f,g_{n}>_{a}|^{2}.
$$
\end{Thm}

{\it Proof.}  Fix $1\le m$ and let
$$
h=\sum_{n=1}^{m}<f,g_{n}>_{a}g_n.
$$
Using the fact that the a-inner product of two functions is
a-periodic (and hence may be factored out of the a-inner product) we
have \bea  <h,h>_{a}  &=& \left <\sum_{n=1}^{m}<f,g_{n}>_{a}g_{n},
\sum_{k=1}^{m}<f,g_{k}>_{a}g_{k}\right >_{a} \no \\ &=&
\sum_{n,k=1}^{m}<f,g_{n}>_{a}
\overline{<f,g_{k}>_{a}}<g_{n},g_{k}>_{a} \no \\
&=&\sum_{n=1}^{m}|<f,g_{n}>_{a}|^{2}. \no  \eea Letting $g=f-h$ we
have by the same type of calculation as above, \bea <h,g>_{a}  &=&
\left <\sum_{n=1}^{m}<f,g_{n}>_{a}g_{n},
f-\sum_{k=1}^{m}<f,g_{k}>_{a}g_{k}\right >_{a} \no \\ &=&
\sum_{n=1}^{m}|<f,g_{n}>_{a}|^{2} -
\sum_{k=1}^{m}|<f,g_{k}>_{a}|^{2} = 0. \no  \eea So we have decomposed
$f$ into two a-orthogonal functions $h,g$.    Therefore,

\bea <f,f>_{a}  &=& <h+g,h+g>_{a}\no \\ &=& <h,h>_{a}+<g,g>_{a} \no \\
 &=& \sum_{n=1}^{m}|<f,g_{n}>_{a}|^{2} + <g,g>_{a} \ge \
 \sum_{n=1}^{m}|<f,g_{n}>_{a}|^{2}. \no  \eea

Since $m$ was arbitrary, we have (2) of the Theorem.  For (1), we just
put together what we know.  By (2) and the Monotone Convergence
Theorem, we have that the series of functions $\sum_{n\in \Bbb
N}|<f,g_{n}>_{a}|^{2}$ converges in $L^{1}[0,a]$.  But, by our
calculations above and the properties of the a-norm,

\bea \|\sum_{n=k}^{m}<f,g_{n}>_{a}g_{n}\|_{L^{2}(\Bbb R)}^{2}
 &=&\int_{0}^{a}\|\sum_{n=k}^{m}<f,g_{n}>_{a}g_{n}\|_{a}^{2}(t)\ dt
 \no \\ &=& \int_{0}^{a}<\sum_{n=k}^{m}<f,g_{n}>_{a}g_{n},
 \sum_{n=k}^{m}<f,g_{n}>_{a}g_{n}>_{a}(t)\   dt \no \\ &=&
 \int_{0}^{a}\sum_{n=k}^{m}|<f,g_{n}>_{a}|^{2}(t)\ dt. \no  \eea Now,
 $\sum_{n\in \Bbb N}|<f,g_{n}>_{a}|^{2}$ converges in $L^{1}[0,a]$
 implies that the right hand side of our equality goes to zero as
 $k\rightarrow \infty$. \\

The ``moreover'' part of the theorem follows immediately from Theorem
\ref{TT} below.  \qed



\section{a-Factorable Operators}\label{Maps}
\setcounter{equation}{0}

Now we consider operators on $L^{2}(\Bbb R)$ which behave naturally
with respect to the a-inner product.  We will call these operators
a-factorable   operators.

\begin{Def}
Fix $E\subset \Bbb R$ and $1\le p\le \infty$.  We say that a linear
operator $L:L^{2}(\Bbb R)\rightarrow L^{p}(E)$ is an {\bf a-factorable
operator} if for any factorization $f={\phi}g$ where $f,g\in
L^{2}(\Bbb R)$   and $\phi$ is an a-periodic function we have
$$
L(f) = L({\phi}g) = {\phi}L(g).
$$
\end{Def}

First we show it is enough to consider factorizations over
$L^\infty([0,a])$

\begin{Prop}
Let $T$ be a bounded operator from $L^{2}(\Bbb R)$ to $L^{2}(E)$.
Then $T$ is a-factorable if and only if $T({\phi}f) = {\phi}T(f)$ for
all $f\in L^{2}(\Bbb R)$ and all a-periodic  ${\phi}\in
L^{\infty}(\mathbb R)$.
\end{Prop}

{\it Proof.}  Assume $\phi$ is a-periodic, $f,g\in L^{2}(\Bbb R)$ and
$f={\phi}g$.  For all $n\in \Bbb N$ let
$$
F_{n} = \{t\in [0,a]:  |{\phi}(t)|>n \}.
$$
Let $E_{n} = [0,1]-F_{n}$ and
$$
\tilde{E}_{n} = \cup_{m\in \Bbb Z}(E_{n}+m)\ \ \text{and} \ \
\tilde{F}_{n} = \cup_{m\in \Bbb Z}(F_{n}+m).
$$
Now,

\bea \|{\chi}_{\tilde{E}_{n}}{\phi}g - {\phi}g\|^{2}_{L^{2}(\Bbb R)}
 &=& \int_{\Bbb R}|{\chi}_{\tilde{F}_{n}}{\phi}(t)g(t)|^{2}\ dt \no \\
 &=& \int_{0}^{a}|{\chi}_{F_{n}}{\phi}(t)|^{2}<g,g>_{a}(t)\ dt. \no
 \eea

Since ${\phi}g\in L^{2}(\Bbb R)$ and $\lim_{n\rightarrow \infty}
{\lambda}(F_{n}) = 0$, it follows that $h_{n}=:
{\chi}_{\tilde{E}_{n}}{\phi}g$ converges to ${\phi}g$ in $L^{2}(\Bbb
R)$.  Since $T$ is a bounded linear operator, it follows that
$T(h_{n})$ converges to $T({\phi}g)$.  But, $T(h_{n}) =
{\chi}_{\tilde{E}_{n}}{\phi}T(g)$ by our assumption.  Now,
$$
\|T(h_{n})\|\le \|T\|\|h_{n}\| \le \|T\|\|{\phi}g\| = \|T\|\|f\|.
$$
Finally, since $|T(h_{n})|\uparrow |{\phi}T(g)|$  it follows from
the Lebesgue Dominated Convergence Theorem that ${\phi}T(g)\in
L^{2}(\Bbb R)$ and $T(h_{n}) \rightarrow {\phi}T(g)$.  This completes
the proof of the   Proposition.  \qed \vspace{14pt}

We have immediately,

\begin{Cor}\label{New}
An operator $T:L^{2}(\Bbb R)\rightarrow L^{p}(E)$ is a-factorable if and
only if $T(E_{m/a}g) = E_{m/a}T(g)$, for all $m\in \Bbb Z$.  That is,
$T$ is a-factorable if and only if it commutes with $E_{m/a}$.
\end{Cor}

The a-inner product naturally defines several types of a-linear maps.
We present two of them here.

\begin{Prop}
Fix $g\in L^{2}(\Bbb R)$ and define a linear operator $L:L^{2}(\Bbb R)
\rightarrow L^{1}[0,1]$ by
$$
L(f) = <f,g>_{a}.
$$
Then $L$ is a bounded, linear a-factorable operator with
$$
\|L\| = \|g\|_{L^{2}(\Bbb R)}.
$$
\end{Prop}

{\it Proof.}  We have that $L$ is a-factorable by Proposition
\ref{a-per}.  Now, for any $f\in L^{2}(\Bbb R)$ we have

\bea \|Lf\| &=& \|<f,g>_{a}\|_{L^{1}[0,a]} \no \\ &=&
 \int_{0}^{a}|\sum_{n\in \Bbb Z} f(t-na)\overline{g(t-na)}|\ dt \no \\
 &\le & \int_{0}^{a}\sqrt{\sum_{n\in \Bbb Z}|f(t-na)|^{2}}
 \sqrt{\sum_{n\in \Bbb Z}|g(t-na)|^{2}}\ dt \no \\ &\le & \left (
 \int_{0}^{a}\sum_{n\in \Bbb Z}|f(t-na)|^{2}\right ) ^{1/2} \left (
 \int_{0}^{a}\sum_{n\in \Bbb Z}|g(t-na)|^{2}\right ) ^{1/2} \no \\
 &=&\|f\|_{L^{2}(\Bbb R)}\|g\|_{L^{2}(\Bbb R)}.\no \eea

Letting $g=f$ we see that $\|L(g)\| = \|g\|$ which, combined with the
above, shows that $\|L\| = \|g\|$.  \qed \vspace{14pt}

Now we define another natural class of a-factorable operators.

\begin{Prop}\label{PP5}
If $g\in L^{\infty}_{a}(\Bbb R)$, the operator
$$
L(f) = <f,g>_{a},
$$
is a bounded linear operator mapping $L^{2}(\Bbb R)$ onto $L^{2}[0,a]$
and
$$
\|L\|^{2} = \text{ess sup}_{[0,a]}<g,g>_{a}.
$$
\end{Prop}

{\it Proof.} This follows directly from the first part of the
proof of Proposition \ref{prebessel} and  again,
 letting $g=f$ above gives the norm of the operator.  \qed
\vspace{14pt}

Now, let $L$ be any a-factorable linear operator from $L^{2}(\Bbb R)$
to $L^{p}(A)$, and let $E=\text{ker}\ L$. If $f\in E$, and $\phi \in
L^{\infty}_{a}(\Bbb R)$, then $L({\phi}f) = {\phi}L(f) = 0$.  So
${\phi}f\in E$.  We summarize this below.

\begin{Prop}\label{PP3}
If $L$ is any a-factorable linear operator with kernel $E$, then $E$
is an a-periodic closed set and so $E^{\perp}= E^{{\perp}_{a}}$.
\end{Prop}

On more property of a-factorable operators into $L^2[0,a]$
 is that the operator  is bounded pointwise by its operator norm with
respect to the a-norm.

\begin{Prop}\label{normlemma}
Let $L$ be an bounded a-factorable linear operator from $L^{2}(\Bbb
R)$ to $L^{2}[0,a]$.  Then for all $f\in L^{2}(\Bbb R)$ we have
$$
|L(f)(t)| \le \|L\|\|f\|_{a}(t),\ \ \text{for all}\ \ t\in [0,a].
$$
\end{Prop}

{\it Proof.}  If not, there is an $f\in L^{2}(\Bbb R)$ and a set
$B \subset [0,a]$ of   positive measure so that

\[
|L(f)(t)|> \|L\|\|f\|_{a}(t),\ \ \text{for all}\ \ t\in B.
\]

{\it  Case 1} If $\|f\|_a(t) =0$ for a.e. $t \in B$. Let
$\Phi=\sum_nT_{na}{\bf 1}_B$ so $\Phi f=0$ yet $L(\Phi f) \ne 0$
 and we have our contradiction.\\

{\it Case 2} If $\|f\|_a(t) \ne 0$ for a.e. $t \in B$ and let
$A \subset B$ so that  $\|f\|_a(t) \ne 0$ for  $t \in A$.  We define 
$\phi=\sum_nT_{na}{\bf 1}_A$.
Now ${\phi}f\in
L^{2}(\Bbb R)$ and
$$
\|\frac{{\phi}f}{\|{\phi}f\|_{a}}\|^2_{L^{2}(\Bbb R)}\le {\lambda}(A).
$$
But,
$$
\|L\left ( \frac{{\phi}f}{\|{\phi}f\|_{a}}\right ) \|^2_{L^2[0,a]}
\ge \int_{A}|L\left ( \frac{{\phi}f}{\|{\phi}f\|_{a}}\right ) (t) |^2
 dt >
{\lambda}(A)\|L\|^2,
$$
which is a contradiction.  \qed \vspace{14pt}

Now we present a short proof of the Riesz representation theorem
for a-factorable operators from $L^2(\mathbb R)$ to $L^1[0,a]$.

\begin{Thm}[Riesz Representation Theorem]
Let $L$ be a bounded a-factorable linear operator from $L^{2}(\Bbb R)$
to $L^{1}[0,a]$.  There exists a function $g\in L^{2}(\Bbb R)$ such
that $L(f) = <f,g>_{a}$, for all $f\in L^{2}(\Bbb R)$ and $\|L\| =
\|g\|_{L^{2} (\Bbb R)}$.
\end{Thm}

{\it Proof.}
Let $f\in L^{2}(\mathbb R)$  and consider the $a$-orthonormal basis
 $ g_n(x) = T_{na}\chi_{[0,a)}(x)$. Hence we have the decomposition
$$
f = \sum_n<f,g_n>_a g_n
$$
We will show the function below is the one we are looking for.
$$
g=\sum_{k\in \mathbb Z}\overline{\widetilde {L(g_k)}}g_k,
$$
where $\widetilde {L(g_k)}$ denotes the periodic extension
of $L(g_k)$ to $\mathbb R$.
First we must show this function is in $L^2(\mathbb R)$.
For positive integers $n$ we define:
$$
h_n=\sum_{|k|\le n} \overline{\widetilde {L(g_k)}}g_k.
$$
For any $\phi \in L^2[0,a]$ we have
$$<\phi,L(g_k)>= \int_0^a \phi (t) L(g_k)(t)\ dt=
\int_0^a L( \tilde \phi g_k)(t)\ dt \le \|\phi\|_{L^2[0,a]} \, \|L\|
$$
Since $\phi$ was arbitrary, $L(g_k) \in L^2[0,a]$ and thus
$\widetilde {L(g_k)}g_k \in L^2(\mathbb R)$.  It follows that
$h_n  \in L^2(\mathbb R)$. Note that
$ \|h_n\|^2_{L^2(\mathbb R)} =\sum_{|k|\le n}
\|L(g_k)\|^2_{L^2[0,a]}$.  Now we compute

\bea
\left \| L \left ( \frac {h_n}
{\|h_n\|_{L^2(\mathbb R)}}\right)\right\|_{L^1[0,a]}
 &=& \frac{1}{\|h_n\|_{L^2(\mathbb R)}}
 \int_0^a\sum_{|k|\le n} |L(g_k)|^2(t)\ dt\no \\
 &=& \|h_n\|_{L^2(\mathbb R)}\le \|L\|. \no \eea
Since $n$ was arbitrary it follows that $g \in L^2(\mathbb R)$.

A direct calculation shows that this is the correct $g$.  i.e.
For all $f\in L^{2}(\mathbb R)$ we have
\bea
 <f,g>_{a}&=& \left<\sum_n<f,g_n>_a g_n,
\sum_k\overline{\widetilde {L(g_k)}}g_k.\right>_a \no \\
 &=&\sum<f,g_n>_a L(g_n)=L(f) \no \eea

\qed \vspace{14pt}

Without much difficulty one may extend this characterization to
a-factorable operators on other $L^p(\mathbb R)$ spaces as well as
into other $L^p[0,a]$ spaces.  We state one of these because it will
be of use in applications to Weyl-Heisenberg frames.

\begin{Prop} Let $L$ be a bounded a-factorable linear operator from
 $L^{2}(\Bbb R)$
to $L^{2}[0,a]$.  There exists a function $g\in L^{2}(\Bbb R)$ such
that $L(f) = <f,g>_{a}$, for all $f\in L^{2}(\Bbb R)$.
\end{Prop}

{\it Proof}.  We note that $L^2[0,a] \subset L^1[0,a]$ and apply
the same proof as above only now it is clear that
$h_n \in L^2(\mathbb R)$.
\qed \vspace{14pt}

We end this section by verifying that for a-factorable operators $T$,
the operator $T^{*}$ behaves as it should relative to the a-inner
product.

\begin{Prop}
If $T$ is an a-factorable operator from $L^{2}(\Bbb R)$ to $L^{2}(\Bbb
R)$, then for all $f,g\in L^{2}(\Bbb R)$ we have

$$
<T^{*}(f),g>_{a} = <f,T(g)>_{a}.
$$

\end{Prop}

{\it Proof.}  Since the a-inner product is a-periodic, we only need to
show the above equality with these functions restricted to
$L^{2}[0,a]$.  For all $m\in \Bbb Z$ we have,

\bea  <T^{*}(f),E_{m/a}g>  &=& \int_{\Bbb R}T^{*}(f)(t)
 \overline{g(t)}e^{-2{\pi}i(m/a)t}\ dt \no \\ &=&
 \int_{0}^{a}<T^{*}(f),g>_{a}(t)e^{-2{\pi}i(m/a)t}\ dt. \no  \eea

Also,

\bea <f,T(E_{m/a}g)> &=& <f,E_{m/a}T(g)> \no \\ &=& \int_{\Bbb
R}f(t)\overline{T(g)(t)}e^{-2{\pi}i(m/a)t}\ dt \no \\ &=&
\int_{0}^{a}<f,T(g)>_{a}e^{-2{\pi}i(m/a)t}\ dt. \no  \eea

Since $<f,T(E_{m/a}g)> = <T^{*}(f),E_{m/a}g>$, for all $m\in \Bbb Z$,
it follows from the above that,
$$
\int_{0}^{a}<T^{*}(f),g>_{a}e^{-2{\pi}i(m/a)t}\ dt =
\int_{0}^{a}<f,T(g)>_{a}e^{-2{\pi}i(m/a)t}\ dt,
$$
for all $m\in \Bbb Z$.  But, this means that
$$
<<T^{*}(f),g>_{a},e^{-2{\pi}i(m/a)t}> =
<<f,T(g)>_{a},e^{-2{\pi}i(m/a)t}>,
$$
for all $m\in \Bbb Z$, where the outer inner product is taken in
$L^{2}[0,a]$.  Since $(\frac{1}{\sqrt{a}}e^{-2{\pi}i(m/a)t})_{m\in  
\Bbb Z}$ is an
orthonormal basis for $L^{2}[0,a]$, we get the desired equality.  \qed
\vspace{14pt}



\section{Weyl-Heisenberg frames and the a-inner product}\label{WH}
\setcounter{equation}{0}

Now we apply our a-inner product theory to Weyl-Heisenberg frames.
For any WH-frame $(g,a,b)$, it is well known that the frame operator
$S$ commutes with $E_{mb}, T_{na}$.  Thus, Corollary \ref{New} yields:

\begin{Cor}
If $(g,a,b)$ is a WH-frame, then the frame operator $S$ is a
1/b-factorable operator.
\end{Cor}

We next
show that the WH-Frame Identity for $(g,a,b)$
 has an interesting  representation in
both the $a$ and the $\frac{1}{b}$ inner products.  The known WH frame
identity requires that the function $f$ be bounded and of compact
support.  While this remains a condition for the WH-Frame Identity
derived from the $a$-inner product we are able to extend this result
to all $f \in L^2( \mathbb R)$ when we use the $\frac{1}{b}$-inner
product.  For this reason we present the theorems separately. \\

The proof of both these theorems have their roots in  the Heil and
Walnut proof of the WH-Frame Identity (see \cite{HW}, Theorem 4.1.5).
We refer   the reader to Proposition \ref{LBD} and Corollary
\ref{sums} for questions concerning convergence of the series and
integrals below.

\begin{Thm}
Let $g\in L^{\infty}_a(\Bbb R)$, and $a,b\in \Bbb R ^+$. For all $f\in
L^{2}(\Bbb R)$ 	which are bounded and compactly supported we have

 $$
\sum_{m,n\in \Bbb Z}|<f,E_{mb}T_{na}g>|^{2} = b^{-1} \sum_k \int_0^a
 <T_{k/b}f,f>_a <g,T_{k/b}g>_a dt.
$$

\end{Thm}

{\it Proof.}  We  start with the WH-frame Identity realizing that
$<g,T_{k/b} g>_a$ is $a$-periodic.
$$
\sum_{m,n\in \Bbb Z}|<f,E_{mb}T_{na}g>|^{2}=
$$

\bea    &=& b^{-1} \sum_k
 \int_{\mathbb R} \overline {f(t)}f(t-k/b) \sum_n g(t-na)
 \overline{g(t-na-k/b)}dt \no \\ &=&  b^{-1} \sum_k \sum_j \int_0^a
 \overline{f(t-ja)} f(t-k/b-ja) <g,T_{k/b}g>_a dt \no \\ &=&  b^{-1}
 \sum_k \int_0^a <T_{k/b}f,f>_a <g,T_{k/b}g>_a dt \no  \eea   \qed
\vspace{14pt}

For the rest of this section we concentrate on the $\frac{1}{b}$ inner
product and its relationship to WH-frames.  In a forthcoming paper  
on the
WH-Frame identity we more closely examine the role of the $a$-inner  
product.
We also show in this paper  that one  may relax the condition on
$g$.  That
 is, the original WH-frame identity holds for all $g \in
L^2(\mathbb R)$  when
 $f$ is bounded and compactly supported.

\begin{Thm}\label{6.2}
Let $g\in L^{\infty}_a(\Bbb R)$, and $a,b\in \Bbb R ^+$. For all $f\in
L^{2}(\Bbb R)$  we have

$$
\sum_{m\in \Bbb Z}|<f,E_{mb}T_{na}g>|^{2} =
\|<f,T_{na}g>_{1/b}\|^{2}_{L^{2}[0,1/b]},
$$
and so
$$
\sum_{n,m\in \Bbb Z}|<f,E_{mb}T_{na}g>|^{2} =
\sum_{n\in
\Bbb Z}\|<f,T_{na}g>_{1/b}\|^{2}_{L^{2}[0,1/b]}.
$$
\end{Thm}

{\it Proof.}
 We just compute,

 \bea
 \sum_{m\in \Bbb Z}|<f,E_{mb}T_{na}g>|^{2} &=&
 \sum_{m\in \Bbb Z}|\int_{\Bbb  R}
 	f(t)\overline{g(t-na)}e^{-2{\pi}imbt}dt|^{2}\no \\
&=& b^{-1}\sum_{m\in \Bbb Z}
	|\sum_{k\in \Bbb Z}\int_{0}^{1/b}f(t-k/b)
	\overline{g(t-na-k/b)}e^{-2{\pi}imbt}dt|^{2} \no \\
&=&b^{-1}\sum_{m\in \Bbb Z}
	|\int_{0}^{1/b}<f,T_{na}g>_{\frac{1}{b}}(t)
	e^{-2{\pi}imbt}dt|^{2} \no \\
&=& b^{-1}\int_{0}^{1/b}
	|<f,T_{na}g>_{\frac{1}{b}}(t)|^{2}\ dt \no \\
&=&  \|<f,T_{na}g>_{1/b}\|^{2}_{L^{2}[0,1/b]}.
\eea
 \qed
\vspace{14pt}

Comparing the equality from Theorem \ref{6.2} above  to the frame
inequalities we have,

\begin{Cor}
Let $g\in$ {\bf PF}, $a,b \in \mathbb R^+$ and define $L:L^{2}(\Bbb
R)\rightarrow   L^{2}(\Bbb R)$ by
$$
L(f) = \sum_{n\in \Bbb
Z}b^{-1}<f,T_{na}g>_{\frac{1}{b}}{\chi}_{[n/b,(n+1)/b)}.
$$
We have
$$
\|L(f)\|^{2} = \sum_{m,n\in \Bbb Z}|<f,E_{mb}T_{na}g>|^{2}.
$$
Hence, $L$ is a bounded linear operator which is an isomorphism if and
only if $(g,a,b)$ is a WH-frame.  Moreover, if $(g,a,b)$ has frame  
bounds
$A,B$, then $\sqrt{A}\|f\|\le \|L(f)\| \le \sqrt{B}\|f\|$, for all $f\in
L^{2}(\Bbb R)$. Hence, $(g,a,b)$ is a normalized tight frame if and  
only if
$L$ is an isometry.
\end{Cor}

Now we want to directly relate our a-inner product to WH-frames.

\begin{Prop}\label{100}
If $g,h\in L^{\infty}_{1/b}(\Bbb R)$, then for all $f\in L^{2}(\Bbb R)$
we have
$$
\sum_{m\in \Bbb Z}<f,E_{mb}g>E_{mb}h = \frac{1}{b}<f,g>_{1/b}h,
$$
where the series converges unconditionally in $L^{2}(\Bbb R)$.  Hence,
$<f,g>_{1/b}g\in \text{span}\ (E_{mb}g)_{m\in \Bbb Z}$.
\end{Prop}

{\it Proof.}  By Proposition \ref{PP5} we know that $<f,g>_{1/b}\in
L^{2}[0,1/b]$.  Next, for any $m\in \Bbb Z$ we have
$$
<f,E_{mb}g> = \int_{\Bbb R}f(t)\
\overline{g(t)}e^{-2{\pi}imbt}\ dt = \int_{0}^{1/b}<f,g>_{1/b}(t)
e^{-2{\pi}imbt}\ dt.
$$
Therefore, if we restrict ourselves to $L^{2}[0,1/b]$ we have
$$
\sum_{m\in \Bbb Z}<f,E_{mb}g>E_{mb} = \sum_{m\in \Bbb Z}\left (
\int_{0}^{1/b}<f,g>_{1/b}e^{-2{\pi}imbt}\right ) e^{2{\pi}imbt} \ dt
$$
$$
= \frac{1}{b}\sum_{m\in \Bbb Z}\left <
<f,g>_{1/b},\sqrt{b}E_{mb}\right > \
\sqrt{b}E_{mb} = \frac{1}{b}<f,g>_{1/b}.
$$
Once we know that we have this convergence in $L^{2}[0,1/b]$, then
redoing the above on $\Bbb R$ with $h$ inserted proves the result and
convergence in $L^{2}(\Bbb R)$.  \qed \vspace{14pt}

There are several interesting consequences of this proposition.
First we recapture the
following result due to de Boor, Devore, and Ron \cite{BDR1}

\begin{Cor}\label{101}
For $g\in L^{2}(\Bbb R)$ and $b\in \Bbb R$, the orthogonal projection
$P$ of $L^{2}(\Bbb R)$ onto span $(E_{mb}g)_{m\in \Bbb Z}$ is
$$
Pf = \frac{1}{\|g\|^{2}_{1/b}}<f,g>_{1/b}\ g  ,
$$
where if $\|g\|_{1/b}(t) = 0$ then $g(t)=0$ so we interpret
$\frac{g(t)}{\|g\|^{2}_{1/b}(t)}=0$. \end{Cor}

{\it Proof.}
By Proposition \ref{OR}, we have that
$(\sqrt{b}E_{mb}\frac{g}{\|g\|_{1/b}})_{m\in \Bbb Z}$ is an orthonormal
sequence in $L^{2}(\Bbb R)$.  Hence, for all $f\in L^{2}(\Bbb R)$
we have
by Proposition \ref{100}
$$
Pf = \sum_{m\in \Bbb Z}<f,\sqrt{b}E_{mb}\frac{g}{\|g\|_{1/b}}>
\sqrt{b}E_{mb}\frac{g}{\|g\|_{1/b}} =
$$
$$
b\sum_{m\in \Bbb
Z}<f,E_{mb}\frac{g}{\|g\|_{1/b}}>E_{mb}\frac{g}{\|g\|_{1/b}}
$$
$$
= <f,\frac{g}{\|g\|_{1/b}>}_{1/b}>\frac{g}{\|g\|_{1/b}}
 =
\frac{1}{\|g\|^{2}_{1/b}}<f,g>_{1/b}\ g .
$$
  \qed \vspace{14pt}

Combining Theorem \ref{Bessel} and Corollary \ref{101} we have:

\begin{Prop}\label{102}
If $(g_{n})_{n\in \Bbb Z}$ is a 1/b-orthonormal sequence in
$L^{2}(\Bbb R)$,
then
$$
P(f) = \sum_{n\in \Bbb Z}<f,g_{n}>_{1/b}g_{n},
$$
is the orthogonal projection of $L^{2}(\Bbb R)$
onto span
$(E_{mb}g_{n})_{n,m\in \Bbb Z}$
 \end{Prop}

In a paper devoted to the study shift-invariant frames and
shift-invariant Riesz bases \cite{RS1}, Ron and Shen develop a
powerful technique called {\bf fiberization}  to decompose the
preframe operator and it's adjoint into a simple collection of
operators called fibers. They then go on to apply this technique to
Weyl-Heisenberg frames by considering the shift invariant space
generated by the countable set $\{E_{mb}g\}_{m \in \mathbb Z}$
\cite{RS2}.  This allows them to produce their amazing result
regarding the
duality principle and Weyl-Heisenberg frames. By using the $1/b$-inner
product we are able to avoid many of the complicated lattice and
dilation  arguments needed for fiberization. In doing so we   produce 
the type of fiberization of the frame operator for a general system that
they have for the self adjoint system. Finally we note that all of our
results have been done on the space side therefore eliminating any
need for taking inverse Fourier transforms to represent the frame
operator.
What results is a simple fiber representation of the WH-frame operator
which we refer to as  a {\bf compression of the frame operator.}

\begin{Thm}
If $(g,a,b)$ is a WH-frame with frame operator $S$, then $S$ has the
form
$$
S(f) = \frac{1}{b}\sum_{n\in \Bbb Z}<f,T_{na}g>_{1/b}T_{na}g
= \frac{1}{b}\sum_{n\in \Bbb Z}P_{n}f\cdot T_{na}\|g\|^{2}_{1/b},
$$
where $P_{n}$ is the orthogonal projection of $L^{2}(\Bbb R)$ onto
span $(E_{mb}T_{na}g)_{m\in \Bbb Z}$ and
 the series converges unconditionally in $L^{2}(\Bbb R)$.
\end{Thm}

{\it Proof.}  If $(g,a,b)$ is a WH-frame then by Proposition \ref{PP6}
we have that $<g,g>_{1/b}\le B$ a.e.  Now, by definition of the frame
operator $S$ we have

\bea S(f) &=& \sum_{m,n\in \Bbb Z}<f,E_{mb}T_{na}g>E_{mb}T_{na}g \no
 \\ &=& \sum_{n\in \Bbb Z}\left ( \sum_{m\in \Bbb Z}
 <f,E_{mb}T_{na}g>E_{mb}T_{na}g \right ) \no \\ &=&
\frac{1}{b}\sum_{n\in \Bbb
 Z}<f,T_{na}g>_{1/b}\ T_{na}g. \no  \eea

An application of Corollary \ref{101} and Theorem \ref{IP} (10)
completes the proof.
\qed \vspace{14pt}

This simple  representation of the frame operator converges
``super-fast''.  That is, we do not have
to compute   any of the modulation parameters to get $S(f)$. While
this  has  immediately become  a useful tool for deriving new properties
regarding the frame operator, because this compression requires us to
pointwise evaluate infinite sums of functions it has obvious
shortcomings in applications.  \\

Theorem \ref{TightWHF} is a classification of the normalized tight  
WH-frames.
We can now restate this in terms of the a-inner products.

\begin{Prop}\label{103}
Let $(g,a,b)$ be a WH-frame.  the following are equivalent:

(1)  $(E_{mb}T_{na}g)_{n,m\in \Bbb Z}$ is a normalized tight
Weyl-Heisenberg
frame.

(2)  $(\frac{1}{\sqrt{b}}T_{n/b}g)_{n\in \Bbb Z}$ is an orthonormal  
sequence
in the a-inner product.

(3)  We have that $g{\perp}_{a} T_{k/b}g$, for all $k\not= 0$ and
$<g,g>_{a} = b$ a.e.
 \end{Prop}

Putting Corollary \ref{102} and Proposition \ref{103} together we have

\begin{Cor}
If $(g,a,b)$ is a normalized tight Weyl-Heisenberg frame,
then
$$
P(f) = \frac{1}{b^{2}}\sum_{n\in \Bbb Z}<f,T_{k/b}g>_{a}T_{k/b}g
$$
is the orthogonal projection of $L^{2}(\Bbb R)$ onto span $(E_{m/a}
T_{k/b}g)_{n,m\in \Bbb Z}$.
\end{Cor}


\section{Two Theorems on WH-Frames}\label{Two}
\setcounter{equation}{0}

In this section we will use the theory developed above to:
(1)  Classify the $g\in L^{2}(\Bbb R)$ for which $(g,a,b)$
is complete when $ab=1$; and (2)  Give an equivalent formulation of the
necessary condition for $(g,a,b)$ to form a WH-frame given
in Theorem \ref{AB}.

First, we need some notation.  If $g\in L^{2}(\Bbb R)$ and
$a>0$ let
$$
X_{g,a} = \overline{\text{span}}\ (E_{ma}g)_{m\in \Bbb Z}.
$$
If $A\subset [0,a]$ and $\phi \in L^{2}(A)$, we write
$\tilde{\phi}$ for the a-periodic extension of $\phi$ to all
of $\Bbb R$.
If $g\in L^{\infty}_{a}(\Bbb R)$, let
$$
g\tilde{L}^{2}[0,a] = \{\tilde{\phi}g: {\phi}\in L^{2}[0,a] \}.
$$

\begin{Lem}\label{Lem5}
Let $E\subset [0,a]$ and $g\in L^{2}(\Bbb R)$ and $A,B>0$.  The
following
are equivalent:

(1)  $A\le \ <g,g>_{a}\ \le B$ a.e. on $E$.

(2)  $A\|\phi \|^{2} \le \|\tilde{\phi}g\|^{2} \le B \|{\phi}\|^{2}$,
for all $\phi \in L^{2}(E)$.
\end{Lem}

{\it Proof.}
If $E\subset [0,a]$, and ${\phi}\in L^{2}(E)$ then
$$
\|\tilde{\phi}\ g\|^{2} = \int_{\Bbb R}|\tilde{\phi}|^{2}
|g|^{2}\ dt = \int_{E}|{\phi}|^{2}<g,g>_{a}\ dt.
$$
Rephrasing this, we have
$$
A\int_{E}|{\phi}|^{2}\ dt \le \int_{E}|{\phi}|^{2}<g,g>_{a}\ dt
\le B\int_{E}|{\phi}|^{2} \ dt,\ \text{for all}\ \phi \in L^{2}[0,a].
$$
The result is immediate from here.
\qed \vspace{14pt}

\begin{Prop}\label{PPP}
Let $g\in L^{\infty}_{a}(\Bbb R)$, $A,B>0$.
  The following are
equivalent:

(1)  $A\le \ <g,g>_{a}\ \le B$ a.e. on the support of $<g,g>_{a}$.

(2)  $A\|\phi \|^{2} \le \|\tilde{\phi}g\|^{2} \le \|\phi \|^{2}$.

(3)  $X_{g,1/a} = g\tilde{L}^{2}[0,a]$.
\end{Prop}

{\it Proof.}
The equivalence of (1) and (2) is Lemma \ref{Lem5}.

$(1)\Rightarrow (2)$:  Let $h\in X_{g,1/a}$.  Choose $h_{n}\in
L^{2}[0,a]$
$$
h_{n} = \sum_{|k|\le n}a_{k}E_{k/a},
$$
so that $\lim_{n\rightarrow \infty}h_{n}g = h$ in $L^{2}(\Bbb R)$.  
Now, for all $m,n\in \Bbb Z$ we have

\bea
 \|h_{n}-h_{m}\|^{2}_{L^{2}[0,a]}
 &=& \int_{0}^{a}|h_{n}(t)-h_{m}(t)|^{2} \ dt \no \\
 &\le& \frac{1}{A}\int_{0}^{a}|h_{n}(t)-h_{m}(t)|^{2}<g,g>_{a}\ dt  
\no \\
 &=& \frac{1}{A}\int_{\Bbb R}|h_{n}(t)g(t)-h_{m}(t)g(t)|^{2}\ dt
 	\rightarrow 0,\ \ \text{as}\ \ n,m\rightarrow \infty. \no
\eea

Therefore, $(h_{n})$ is Cauchy in $L^{2}[0,a]$ and hence convergent
to some $f\in L^{2}[0,a]$.  Now,

\bea
\|h_{n}g-fg\|_{L^{2}(\Bbb R)}
 &=& \int_{\Bbb	R}|h_{n}(t)g(t) -f(t)g(t)|\ dt \no \\
 &=& \int_{0}^{a}|h_{n}(t)-f(t)|<g,g>_{a}\ dt \no  \\
 &\le&  B\int_{0}^{a}|h(t)-f(t)|^{2}\ dt\rightarrow 0,\ \ \text{as}\ \
	n\rightarrow \infty. \no
\eea

That is, $fg = h\in g\tilde{L}[0,a]$.

$(2)\Rightarrow (1)$:  Define,
\[
g_{n}(t) = \left \{
\begin{array}{cl}
 g(t) : & |g(t)| \le n \\
0 : & \text{otherwise}.
\end{array} \right.
\]

Let
$$
E = (\text{supp}\ <g,g>_{a})\cap [0,a].
$$
Define $T_{n}:L^{2}(E)\rightarrow L^{2}(\Bbb R)$ by $T_{n}(\phi ) =
\tilde{\phi}g_{n}$.  For all $n$, $T_{n}$ is a bounded linear operator
and
$$
\|T_{n}(\phi )\| = \|\tilde{\phi}g_{n}\|\le \|\tilde{\phi}g\|,
$$
Therefore, the $(T_{n})$ are pointwise bounded.  By the Uniform
Boundedness
Principle, the $(T_{n})$ are uniformly bounded.  Also,
$$
|T_{n}\phi | \uparrow |\tilde{\phi} g |\in L^{2}(\Bbb R),
$$
since $g\tilde{L}^{2}[0,a] = X_{g,1/a}\subset L^{2}(\Bbb R)$.
Hence, the operator $T$ defined by $T(\phi )= \tilde{\phi}g$ is a
bounded
linear operator from $L^{2}(E)$ to $g\tilde{L}[0,a]$, which is
one-to-one.
Since
$X_{g,1/a}$ is a Banach space, it follows that $T$ is an
isomorphism.  Hence,
there are constants $A,B>0$ satisfying,
$$
A\|\phi \|^{2} \le \|\tilde{\phi}g\|^{2} \le B\|\phi \|^{2}.
$$
That is,
$$
A\int_{E}|\phi (t)|^{2}\ dt \le \int_{E}|\phi (t)|^{2}(t)<g,g>_{a}\ dt
\le B\int_{E}|\phi (t) |^{2}\ dt.
$$
Hence, $A\le \ <g,g>_{a}\ \le B$.
 \qed \vspace{14pt}

Now we have an important consequence of these results which is a
the modulation version of the shift-invariant result of Ron and Shen
\cite{RS1} (Theorem 2.2.14) .

\begin{Thm}\label{TT}
For $g\in L^{2}(\Bbb R)$, the following are equivalent:

(1)  There are numbers $A,B>0$ so that $A\le \ <g,g>_{a}\ \le B$ a.e.

(2)  $(E_{m/a}g)_{m\in \Bbb Z}$ is a Riesz basic sequence.

Hence, if $ab=1$, then (1) and (2) are equivalent to

(3)  $(E_{mb}g)_{m\in \Bbb Z}$ is a Riesz basic sequence with Riesz
basis constants $\sqrt{A},\sqrt{B}$.
\end{Thm}

{\it Proof.}
$(1)\Rightarrow (2)$:  In the proof of Proposition \ref{PPP} we saw  
that the
map $T:L^{2}[0,a]\rightarrow L^{2}(\Bbb R)$ given by $T(\phi ) =
\tilde{\phi}g$ is an isomorphism.  Since $(\frac{1}{\sqrt{a}}
E_{m/a})_{m\in
\Bbb Z}$ is an orthonormal basis for $L^{2}[0,a]$, it follows that
$(T(E_{m/a})) = (E_{m/a}g)$ is a Riesz basic sequence.

$(2)\Rightarrow (1)$:  By assumption there are constants $A,B>0$
satisfying
for all sequences of scalars $(a_{m})_{m\in \Bbb Z}$,
$$
A\sum_{m\in \Bbb Z}|a_{m}|^{2} \le \|\sum_{m\in \Bbb
Z}|a_{m}|E_{mb}g\|^{2}_{
L^{2}(\Bbb R)}\le B \sum_{m\in \Bbb Z}|a_{m}|^{2}.
$$
Since $\|\sum_{m}|a_{m}|E_{mb}\|^{2} = \sum_{m}|a_{m}|^{2}$, it
follows that
for all $\phi \in L^{2}[0,a]$, $A\|\phi \|^{2} \le
\|\tilde{\phi}g\|^{2}\le
B\|\phi \|$.
 \qed \vspace{14pt}

Note that Theorem \ref{TT} is really half of Theorem \ref{dual}.  This
seems to indicate that there is ``another half'' someplace which
produces
the whole result.  It would be interesting to find this. Note also that
if $g(t) = e^{-t^{2}}$, then the Fourier transform of $g$ is
$\hat{g}(t) = \sqrt{2\pi }\ e^{-t^{2}/2}$.  A direct calculation shows
that there are constants $A,B>0$ such that
$$
A\le \ <g,g>_{1}\ \le B \ \ \text{and}\ \
A\le \ <\hat{g},\hat{g}>_{1}\ \le B .
$$
It follows that $(E_{m}\hat{g}) = (\hat{T_{m}g})$ is a Riesz basic
sequence.  So it follows that both $(E_{m}g)_{m\in \Bbb Z}$ and
$(T_{m}g)_{m\in \Bbb Z}$ are Riesz basic sequences, despite the fact
that $(g,a,b)$ is not a WH-frame.\\

Next we will completely identify the functions $g\in L^{2}(\Bbb R)$ for
which $(g,a,b)$ is complete in $L^{2}(\Bbb R)$.  We need one more
piece of
notation.  If $(E_{n})_{n\in \Bbb Z}$ is any orthonormal basis for
$L^{2}[0,a]$, we let $R$ denote the {\bf right hand shift operator}  
given
by $R(E_{n}) = E_{n+1}$, for all $n\in \Bbb Z$.

\begin{Prop}\label{Late}
Let $E_{n} = e^{2{\pi}int} \in L^{2}[0,1]$, for all $n\in \Bbb Z$, and
let $f = \sum_{n\in \Bbb Z}a_{n}E_{n}\in L^{2}[0,1]$.  The following are
equivalent:

(1)  $(R^{n}f)_{n\in \Bbb Z}$ is complete in $L^{2}[0,1]$.

(2)  $f\not= 0$, a.e. in $L^{2}[0,1]$.
\end{Prop}

{\it Proof.}
First we compute,
$$
R(f) = R\left ( \sum_{n\in \Bbb Z}a_{n}E_{n}\right ) = \sum_{n\in
\Bbb Z}
a_{n}E_{n+1} = E_{1}\sum_{n\in \Bbb Z}a_{n}E_{n} = E_{1}f.
$$
Note that for $h\in L^{2}[0,1]$, we have that $h\perp E_{n}f$ if
and only if
$h\overline{f} \perp E_{n}$.  Hence, $h\perp E_{n}f$, for all $n\in  
\Bbb Z$
if and only if $h\overline{f} = 0$, a.e.  It follows that
$(E_{n}f)_{n\in \Bbb Z}$ is complete (and hence $(R^{n}f)$ is
complete) if
and only if we have:  Whenever $h\in L^{2}[0,1]$ and $h\overline{f}  
= 0$,
a.e., then $h=0$ a.e.  This is clearly equivalent to $f\not= 0$, a.e.
 \qed \vspace{14pt}

Now we can give the required classification.  If $f(x,y)$ is a function
of two variables, we write $f_{x}$ for the function $f_{x}(y) = f(x,y)$
and $f_{y}$ for the function $f_{y}(x) = f(x,y)$.

\begin{Thm}\label{Main}
Let $a=b=1$ and $g\in L^{2}(\Bbb R)$.  The following are equivalent:

(1)  $(E_{m}T_{n}g)_{m,n\in \Bbb Z}$ is complete in $L^{2}(\Bbb R)$.

(2)  There is a function $f(x,y):[0,1]\times [0,1]
\rightarrow \Bbb R$ satisfying:

\ \ \ \ (a)  For a.e. $y\in [0,1]$, we have that $f_{y}\not= 0$, a.e.

\ \ \ \ (b)  For all $y\in [0,1]$, we have
$$
g(y+n) = \hat{f}_{y}(n),\ \ \text{for all}\ \ n\in \Bbb Z ,
$$
where $\hat{f}_{y}(n) = <f_{y},E_{n}>$.
\end{Thm}

{\it Proof.}
$(2)\Rightarrow (1)$:  Suppose that $h\in L^{2}(\Bbb R)$ and
$h\perp \text{span}\ (E_{mb}T_{na}g)_{n,m\in \Bbb Z}$.  Then by
Proposition \ref{PP1},
$$
<h,T_{na}g>_{1/b} = <f,T_{na}>_{a} = 0,\ \text{a.e. for all n.}
$$
That is,
$$
\sum_{k\in \Bbb Z}f(y-ka)\overline{g(y-(k-n)a)} = 0,\ \text{a.e.
for all n.}
$$
Letting $h_{y} = \sum_{k}h(y-ka)e^{2{\pi}ix}$ and $g_{y} = \sum_{k}
g(y-ka)e^{2{\pi}ix}$, we have that $g_{y}(x) = f_{y}(x)$.  Also by the
above we have that
$$
h_{y}\perp R^{n}(f_{y}),\ \text{for all n} \ \in \Bbb Z.
$$
Hence, by Proposition \ref{Late}, we have that $h_{y} = 0$ a.e.
That is,
$h(y-ka) = 0$, for all $k\in \Bbb Z$.  Hence, $h=0$ a.e. and it follows
that $(g,a,b)$ is complete.

$(1)\Rightarrow (2)$:  Define the function
$$
f(x,y) = \sum_{k\in \Bbb Z}g(y-ka)e^{2{\pi}ikx}.
$$
Then the above argument for $(2)\Rightarrow (1)$ shows

 that $f(x,y)$
has the desired properties.
   \qed \vspace{14pt}

Recall \cite{KZ} that a class of infinitely differentiable functions on 
$\Bbb T$ is called {\bf quasi-analytic} if the only function in the  
class which
vanishes with all its derivatives at some point $t_{0}\in \Bbb T$ is the
function which vanishes identically.  A direct calculation shows
that functions
in a quasi-analytic class can have at most a finite number of zeroes 
on $\Bbb T$.  On page 113 of Katznelson \cite{KZ} is the following  
theorem.

\begin{Thm}(Denjoy-Carleman)
Given some $d,K>0$ in $\Bbb R$, let
$$
E = \{f= \sum_{n\in \Bbb N}a_{n}e^{2{\pi}int} : |a_{n}|\le Kd^{n}, \ 
\text{for all}\ n\}.
$$
Then $E$ is a quasi-analytic class.
\end{Thm}

Combining the above we have

\begin{Thm}
Let $g\in L^{2}(\Bbb R)$ and assume there exist $K,d>0$ so that
$$
|g(t+n)|\le Kd^{n},\ \text{for all}\ t\in [0,a].
$$
Then $(E_{m}T_{n})_{m,n\in \Bbb Z}$ is complete in $L^{2}(\Bbb R)$.
In particular, $g(t) = e^{-ct^{2}}$, works for all $c>0$.
\end{Thm}

Finally, we recall \cite{Duren} that if $f$ is an $H^{p}$-function
then $log|f(e^{i{\theta}})|$ is integrable unless $f(z)\equiv 0$.
In particular, if $f$ vanishes on a set of positive measure then it
vanishes identically.  One consequence of this is that if $m\in \Bbb Z$
and
$f = \sum_{k=m}^{\infty} a_{k}E_{k}\in L^{p}[0,1]$
and $a_{i}\not= 0$ for at least one $m\le i < \infty$, then $f\not=  
0$ a.e.
Combining this with the proof of Theorem \ref{Main} (the proof of
$(2)\Rightarrow (1)$) we have

\begin{Thm}
Let $g\in L^{2}(\Bbb R)$ be supported on a ray $[{\alpha},\infty)$
(In particular, if $g$ has compact support).
The following are equivalent:

(1)  $(E_{m}T_{n}g)_{m,n\in \Bbb Z}$ is complete in $L^{2}(\Bbb R)$.

(2)  $\text{sup}_{n\in \Bbb Z}|g(x-n)| \not= 0$ a.e.
\end{Thm}



\section{Frames in the a-inner product}\label{a-frames}
\setcounter{equation}{0}

We will now look at the notion of frames and Riesz bases in the
a-inner product.

\begin{Def}
We say that a sequence $f_{n}\in L^{2}(\Bbb R)$ is a {\bf
a-Riesz basic sequence} if there is an a-orthonormal basis
$(g_{n})_{n\in \Bbb Z}$ and an a-factorable operator $T$ on
$L^{2}(\Bbb R)$ with $T(g_{n}) = f_{n}$ so that $T$ is invertible
on its range.  If $T$ is surjective, we call $(f_{n})$
a {\bf a-Riesz basis} for $L^{2}(\Bbb R)$.
\end{Def}

\begin{Prop}
For $f_{n}\in L^{2}(\Bbb R)$, for all $n\in \Bbb Z$, the following
are equivalent:

(1)  $(f_{n})_{n\in \Bbb Z}$ is an a-Riesz basic sequence.

(2)  $(E_{m/a}f_{n})_{n\in \Bbb Z}$ is a Riesz basic sequence.
\end{Prop}

{\it Proof.}
$(1)\Rightarrow (2)$:  By assumption, there is an a-orthonormal
basis $(g_{n})$ and an a-factorable operator $T$ with
$$
T(g_{n}) = f_{n},\ \ \text{for all} \ \ n\in \Bbb Z .
$$
By the definition of an a-orthonormal basis we have that
$(\frac{1}{\sqrt{a}}E_{m/a}g_{n})_{m,n\in \Bbb Z}$ is an
orthonormal basis $L^{2}(\Bbb R)$.  Since $T$ is an isomorphism,
it follows that
$$
(T(\frac{1}{\sqrt{a}}E_{m/a}g_{n}))_{n,m\in \Bbb Z}
 =
(\frac{1}{\sqrt{a}}E_{m/a}T(g_{n}))_{n,m\in \Bbb Z}
= (\frac{1}{\sqrt{a}}E_{m/a}f_{n})_{n,m\in \Bbb Z}
$$
is a Riesz basic sequence.

$(2)\Rightarrow (1)$:  Let $g={\chi}_{[0,a)}$ so that
$(\frac{1}{\sqrt{a}}E_{m/a}T_{na}g)_{m,n\in \Bbb Z}$ is an
orthonormal basis for $L^{2}(\Bbb R)$.  Then
$$
T(\frac{1}{\sqrt{a}}E_{m/a}T_{na}g)= E_{m/a}f_{n}
$$
is an a-factorable linear operator which is an isomorphism
because $(E_{m/a}f_{n})$ is a Riesz basic sequence.  Hence,
$(f_{n})$ is an a-Riesz basic sequence.
\qed \vspace{14pt}

\begin{Cor}
For $f_{n}\in L^{2}(\Bbb R)$, for all $n\in \Bbb Z$, the following
are equivalent:

(1)  $(f_{n})_{n\in \Bbb Z}$ is an a-Riesz basis.

(2)  $(E_{m/a}f_{n})_{n\in \Bbb Z}$ is a Riesz basis for
$L^{2}(\Bbb R)$.
\end{Cor}

Since the inner product on a Hilbert space is used to define a frame,
we can get a corresponding concept for the a-inner product.

\begin{Def}
If $g_{n}\in L^{2}(\Bbb R)$, for all $n\in \Bbb Z$, we call
$(g_{n})_{n\in \Bbb Z}$ an {\bf a-frame sequence} if there exist
constants $A,B>0$ so that for all
$f\in \text{span}\ (E_{m/a}g_{n})_{m,n\in
\Bbb Z}$ we have
$$
A\|f\|_{a}^{2}\le \sum_{n\in \Bbb Z}|<f,g_{n}>_{a}|^{2}\le
B\|f\|_{a}^{2}.
$$

If the inequality above holds for all $f \in L^2(\mathbb R)$ then
we call $g_n$ an  {\bf $a$-frame}.
\end{Def}

Now we have the corresponding result to Theorem \ref{T1}.

\begin{Thm}\label{a-frame}
Let $g_{n}\in L^{2}(\Bbb R)$, for all $n\in \Bbb Z$.
The following are equivalent:

(1)  $(g_{n})_{n\in \Bbb Z}$ is an a-frame.

(2)  If $(e_{n})_{n\in \Bbb Z}$ is an a-orthonormal basis for
$L^{2}(\Bbb R)$, and $T:L^{2}(\mathbb R)\rightarrow L^{2}(\mathbb R)$
with
$T(e_{n}) = g_{n}$ is a-factorable, then
$T$ is a bounded, linear surjective operator on $L^{2}(\Bbb R)$.
\end{Thm}

{\it Proof.}
If $T(e_{n}) = g_{n}$, then
$$
<T^{*}(f),e_{n}>_{a} = <f,T(e_{n})>_{a} = <f,g_{n}>_{a}.
$$
Hence, by Theorem \ref{Bessel} we have that
$T^{*}(f) = \sum_{n\in \Bbb Z}<f,g_{n}>_{a}e_{n}$ and
$$
\|T^{*}(f)\|_{a}^{2} = \sum_{n\in \Bbb Z}|<f,g_{n}>_{a}|^{2}.
$$
Hence, $(g_{n})$ is an a-frame sequence if and only if
$$
A\|f\|_{a}^{2} \le \|T^{*}(f)\|_{a}^{2}\le B\|f\|_{a}^{2},\ \
\text{for all}\
\ f\in L^{2}(\mathbb R).
$$
But this is equivalent to $T^{*}$ being an isomorphism, which
itself is equivalent to $T$ being a bounded, linear onto operator.
\qed \vspace{14pt}

Finally, we can relate this back to our regular frame sequences.

\begin{Prop}\label{a-frame sequence}
Let $g_{n}\in L^{2}(\Bbb R)$, for all $n\in \Bbb Z$.  The following
are equivalent:

(1)  $(g_{n})_{n\in \Bbb Z}$ is an a-frame sequence.

(2)  $(E_{m/a}g_{n})_{m,n\in \Bbb Z}$ is a frame sequence.
\end{Prop}

{\it Proof.}
$(1)\Rightarrow (2)$:  If $(g_{n})$ is an a-frame sequence, then
there is
an a-orthonormal basis $(e_{n})$ for $L^{2}(\Bbb R)$ and an a-factorable
onto (closed range) operator $T(e_{n}) = g_{n}$.  Now,
$(E_{m/a}e_{n})_{n,m\in \Bbb Z}$ is an orthonormal basis for
$L^{2}(\Bbb R)$
and $$
T(E_{m/a}e_{n}) = E_{m/a}T(e_{n}) = E_{m/a}g_{n}.
$$
Hence, $(E_{m/a}g_{n})_{m,n\in \Bbb Z}$ is a frame sequence.

$(2)\Rightarrow (1)$:  Reverse the steps in part I above.
\qed \vspace{14pt}

The following Corollary is immediate from Theorem \ref{a-frame} and
Proposition \ref{a-frame sequence}.

\begin{Cor}
Let $g\in L^{2}(\mathbb R)$ and $a,b\in \mathbb R$.  The following are
equivalent:

(1)  $(g,a)$ is a 1/b-frame.

(2)  $(g,a,b)$ is a Weyl-Heisenberg frame.
\end{Cor}



\section{Gram-Schmidt Process}\label{GS}
\setcounter{equation}{0}

In this section we will look at the Gram-Schmidt process for the
a-inner product.  First we need a result which shows that this
process produces functions which are in the proper spans.

\begin{Prop}\label{P20}
Let $f,g,h\in L^{2}(\Bbb R)$.  We have:

(1)  $N_{a}(g)\in \text{span}\ (E_{m/a}g)_{m\in \Bbb Z}$.

(2)  If any two of $f,g,h$ are
in $L^{\infty}_{a}(\Bbb R)$, then $<f,h>_{a}g\in \text{span}\
(E_{m/a}g)_{m\in
\Bbb Z}$.
\end{Prop}

{\it Proof.}

$(1)$:  For each $n\in \Bbb N$ let
$$
E_{n} = \{t\in [0,a]:|<g,g>_{a}(t)|^{2}\ge n\ \text{or}\
<g,g>_{a}(t) \le \frac{1}{n} \}.
$$
Also, let
$$
\tilde{E}_{n} = \cup_{m\in \Bbb Z}(E_{n}+m).
$$
Since $g\in L^{2}(\Bbb R)$, we have
$$
\|g\|^{2} = \int_{0}^{a}<g,g>_{a}(t)\ dt < \infty.
$$
Hence, $\lim_{n\rightarrow \infty}{\lambda}(E_{n}) = 0$.  Let
$F_{n} = [0,a] - E_{n}$ and
$$
\tilde{F}_{n} = \cup_{m\in \Bbb Z}(F_{n}+m).
$$
Now,
$$
\frac{1}{n}\le
\frac{1}{<{\chi}_{\tilde{F}_{n}}g,{\chi}_{\tilde{F}_{n}}g>_{a}}
\le n.
$$
Hence,
$$
\frac{1}{<{\chi}_{\tilde{F}_{n}}g,{\chi}_{\tilde{F}_{n}}g>_{a}} \in
L^{\infty}_{a}(\Bbb R). $$
Hence,
$$
\frac{{\chi}_{\tilde{F}_{n}}g}{<{\chi}_{\tilde{F}_{n}}g,
{\chi}_{\tilde{F}_{n}}g>_{a}}+ {\chi}_{\tilde{E}_{n}}g \in \text{span}\
(E_{m/a}g)_{m\in \Bbb Z}.
$$
Also,

\bea
 \| \frac{{\chi}_{\tilde{F}_{n}}g}{<{\chi}_{\tilde{F}_{n}}g,
 	{\chi}_{\tilde{F}_{n}}g>_{a}} + {\chi}_{\tilde{E}_{n}}g
 	-N_{a}(g)\|_{L^{2}(\Bbb R)}
 &=& \| {\chi}_{\tilde{E}_{n}}g -
	\frac{{\chi}_{\tilde{E}_{n}}g }{<g,g>_{a}}\| \no \\
 &\le& \| {\chi}_{\tilde{E}_{n}}g\| +
	\|\frac{{\chi}_{\tilde{E}_{n}}g  }{<g,g>_{a}}\| \no \\
 &=& \left ( \int_{\Bbb R}|{\chi}_{\tilde{E}_{n}}g|^{2}\ dt \right  
) ^{1/2}
	+ \|N_{a}({\chi}_{\tilde{E}_{n}}g )\| \no \\
 &\le&  \left ( \int_{E_{n}}<g,g>_{a}(t)\ dt \right ) ^{1/2} +
	{\lambda}(E_{n}). \no
\eea

But the right hand side of the above inequality goes to zero as
$n\rightarrow
\infty$.

$(2)$:  Assume first that $f,h\in L^{\infty}_{a}(\Bbb R)$.  Let
$B=\|f\|+{a}$
and $C = \|h\|_{a}$.  Now,

\bea
 |<f,h>_{a}|
 &=& |\sum_{n\in \Bbb Z}f(t-na)\overline{g(t-na)} \no \\
 &\le& \sqrt{\sum_{n\in \Bbb Z}|f(t-na)|^{2}}
	\sqrt{\sum_{n\in \Bbb Z}|g(t-na)|^{2}} \le
\sqrt{B}\sqrt{C}. \no
\eea

Therefore, $<f,h>_{a}$ is a bounded a-periodic function on $\Bbb R$.
this implies that $<f,h>_{a}g\in L^{2}(\Bbb R)$.

Now suppose that $g,h\in L^{\infty}_{a}(\Bbb R)$.  Let $B=
\|g\|_{a}$ and
$C=\|h\|_{a}$.  Then

\bea
 \|<f,h>_{a}g\|^{2}_{L^{2}(\Bbb R)}
 &=& \|\left ( \sum_{n\in \Bbb Z}f(t-na)
	\overline{h(t-na)}\right ) g\|_{L^{2}(\Bbb R)} \no \\
 &=&\int_{0}^{a}|\sum_{n\in \Bbb Z}f(t-na)\overline{h(t-na)}|^{2}
	\sum_{n\in \Bbb Z}|g(t-na)|^{2}\ dt \no \\
 &\le& B\int_{0}^{a}\sum_{n\in \Bbb Z}|f(t-na)|^{2}
	\sum_{n\in \Bbb Z}|h(t-na)|^{2} \no \\
 &\le&  BC\|f\|^{2}_{L^{2}(\Bbb R)}. \no
\eea

Recall that
$$
\text{span}\ (E_{m/a}g)_{m\in \Bbb Z} = \{{\phi}g:{\phi}\ \text{is  
a-periodic
and}
\ \ {\phi}g\in L^{2}(\Bbb R)\}.
$$
So by the above, we have that $<f,h>_{a}g\in \text{span}\
(E_{m/a}g)_{m\in
\Bbb Z}$.
\qed \vspace{14pt}

\begin{Def}
Let $g_{n}\in L^{2}(\Bbb R)$, for $1\le n\le k$.  We say that
$(g_{n})_{n=1}^{k}$ is {\bf a-linearly independent} if for each
$1\le n\le k$, $g_{n}\notin \text{span}\ (E_{m/a}g_{i})_{m\in \Bbb Z;
1\le i\not= n \le k}$.  An arbitrary family is {\bf a-linearly
independent} if every finite sub-family is a-linearly independent.
\end{Def}

Now we carry out the Gram-Schmidt process.

\begin{Thm}(Gram-Schmidt orthonormalization procedure)
Let $(g_{n})_{n\in \Bbb
N}$ be an a-linearly independent
 sequence in $L^{2}(\Bbb R)$ for $a>0$.  Then there exists
an a-orthonormal sequence $(e_{n})_{n\in \Bbb N}$ satisfying
for all $n\in \Bbb N$:
$$
\text{span}\ (E_{m/a}g_{k})_{m\in \Bbb Z,1\le k \le n} =
\text{span}\ (E_{m/a}e_{k})_{m\in \Bbb Z,1\le k \le n} .
$$
\end{Thm}

{\it Proof}
We proceed by induction.  First let $e_{1} = N_{a}(g_{1})$.
If $(e_{i})_{i=1}^{n}$ have been defined to satisfy the
theorem, let
$$
e_{n+1} = N_{a}(g_{n+1}-\sum_{i=1}^{n}<g_{i},e_{i}>_{a}e_{i}).
$$
Let
$$
h = g_{n+1}-\sum_{i=1}^{n}<g_{n+1},e_{i}>_{a}e_{i}.
$$
Note that $h\not= 0$ by our a-linearly independent assumption and
Proposition \ref{P20}.
Now, for $1\le k \le n$ we have

\bea
 <e_{n+1},e_{k}>_{a}
 & =& \frac{1}{<h,h>_{a}}\left ( <g_{n+1},e_{k}>_{a}
	- \sum_{i=1}^{n}<g_{n+1},e_{i}>_{a}<e_{i},e_{k}>_{a}\right  
) \no \\
 &=& \frac{1}{<h,h>_{a}}\left ( <g_{n+1},e_{k}>_{a} -
	<g_{n+1},e_{k}>_{a}<e_{k},e_{k}>_{a}\right ) = 0. \no
\eea

The statement about the linear spans follows from Proposition
\ref{P20}.
\qed \vspace{14pt}



\end{document}